\documentclass[3p]{elsarticle}

\usepackage[english]{babel}

\usepackage{amsmath,amssymb}

\usepackage{todonotes}
\newtheorem{assumption}{Assumption}

\usepackage{hyperref}
\usepackage{color}

\usepackage{tikz}
\usetikzlibrary{arrows,shapes,snakes,calendar,matrix,backgrounds,folding,spy,positioning}

\DeclareMathOperator*{\argmin}{\arg\!\min}

\newcommand{\phiexact}{\phi^{\text{ex}}}
\newcommand{\rr}{\mathbb{R}}
\newcommand{\jumpleft}{[\![}
\newcommand{\jumpright}{]\!]}

\newcommand{\averageleft}{\{\!\!\{}
\newcommand{\averageright}{\}\!\!\}}
\renewcommand{\div}{\textrm{div}\ \!}
\newcommand{\mT}{\mathbb{T}} 

\newtheorem{remark}{Remark}[section]

\numberwithin{equation}{section}
\numberwithin{figure}{section}

\journal{Comp. Meth. Appl. Mech. Eng.}

\begin{document}


\begin{frontmatter}
\title{High order unfitted finite element methods on level set domains using isoparametric mappings}
\author{Christoph Lehrenfeld}
\address{Institute for Computational and Applied Mathematics, WWU M\"unster, M\"unster, Germany}

\begin{abstract}
We introduce a new class of unfitted finite element methods with high order accurate numerical integration over curved surfaces and volumes which are only implicitly defined by level set functions. 
An unfitted finite element method which is suitable for the case of piecewise planar interfaces is combined with a parametric mapping of the underlying mesh resulting in an isoparametric unfitted finite element method. The parametric mapping is constructed in a way such that the quality of the piecewise planar interface reconstruction is significantly improved allowing for high order accurate computations of (unfitted) domain and surface integrals. This approach is new. We present the method, discuss implementational aspects and present numerical examples which demonstrate the quality and potential of this method. 
\end{abstract}

\begin{keyword}
numerical integration \sep
level set \sep
unfitted finite element method \sep
isoparametric finite element method \sep
high order methods \sep
interface problems
\end{keyword}


\end{frontmatter}
\section{Introduction}\label{sec:intro}
\subsection{Motivation}\label{sec:intro:motivation}
In the recent years \emph{unfitted finite element methods} have drawn more and more attention. 
 The possibility to handle complex geometries without the need for complex and time consuming mesh generation is very appealing. The methodology of unfitted finite element methods, i.e. methods which are able to cope with interfaces or boundaries which are not aligned with the grid, has been investigated for different problems. For boundary value problems with non-aligned boundaries methods such as penalty methods \cite{babuska73b,barrettelliott86}, the fictitious domain method \cite{glowinskietal94,burman2012fictitious}, the immersed boundary method \cite{peskinmcqueen89}, and other unfitted finite element methods \cite{bastian2009unfitted} have been developed. 
For \emph{unfitted} interface problems extended finite element methods (XFEM) have been developed in (among others) \cite{fries2010extended, hansbo2002unfitted, gross04, massjung12, Becker20093352}. Also partial differential equations on surfaces which are not meshed have been considered using the trace finite element method (TraceFEM) \cite{olshanskii2009finite}. 

In the community of unfitted finite element methods mostly piecewise linear discretizations are considered. One major issue in the design and realization of high order methods is the problem of numerical integration on domains which are only prescribed implicitly, for instance as the zero level of a scalar function, the so-called level set function. The use of standard integration rules (which ignore the cut position on cut elements) is not a good option as the integrand does not provide the necessary smoothness.

The purpose of this paper is the presentation of a new approach which allows for high order numerical integration on domains prescribed by level set functions. The approach is new, robust and fairly simple to implement. The method is geometry-based and can be applied to unfitted interface or boundary value problems as well as to partial differential equations on surfaces.  

\subsection{Literature}\label{sec:intro:lit}
We briefly discuss the state of the art in the literature to put the new approach into context. One approach that is often used in order to realize numerical integration on implicit domains consists out of two step: the approximation of the interface with a piecewise planar interface and a tesselation algorithm to divide the subdomains of a cut element into simple geometries, on which standard quadrature rules can be applied. A famous example for quadrilaterals and cubes is the marching cube algorithm \cite{lorensen1987marching}. For simplices a detailed discussion of this approach can be found in (among others) \cite[Chapter 5]{naerland2014geometrychap5},\cite{mayer2009interface} for triangles and tetrahedra and in \cite{lehrenfeld2015nitsche,lehrenfeld2015diss} also for 4-prisms and pentatopes (4-simplices). Many simulation codes which apply unfitted finite element discretizations, e.g. \cite{DROPS,engwer2012dune,burman2014cutfem,renard2014getfem++,carraro2015implementation} make use of such strategies. In order to increase the accuracy of this integration one often combines this approach with subdivisions and adaptivity. Especially on octree-based meshes this can be done very efficiently \cite{chernyshenko2015adaptive}. However, this tesselation approach is only second order accurate (w.r.t. the finest subtriangulation) which complicates the realization of high order methods. 

One approach to solve this problem is the tailoring of quadrature points and weights which provide high order accurate integration rules for implicit domains.
Such a construction of points and weights is based on the idea of fitting integral moments, cf. \cite{muller2013highly,sudhakar2013quadrature}. Although this results in accurate integration rules it has two shortcomings. First, the computation of quadrature rules is fairly involved. This aspect is typically outwayed by the resulting high accuracy. Secondly, the major problem, quadrature weights are in general not positive. This can lead to stability problems as positivity of mass or stiffness matrices in finite element formulations can no longer be guaranteed.

For special cases satisfactory answers to the question of high order numerical integration strategies have been found which allow for an implementation of high order unfitted finite element methods. We mention a few interesting approaches. 
For quadrilaterals and hexahedra in \cite{saye2015hoquad} a numerical integration algorithm is presented which can achieve arbitrary high order accuracy and guarantee positivity of integration weights at the same time. The approach is based on the idea of interpreting the interface as a piecewise graph over a hyperplanes. In \cite{burman2015cut} an unfitted boundary value problem is considered. Instead of aiming at a high order accurate geometrical approximation of the boundary a correction in the imposition of boundary values is applied in order to recover a high order method. 
For the discretization of partial differential equations on surfaces, in \cite{grande2014highorder} 
a \emph{parametric mapping of the interface} from a piecewise planar representation to a high order representation based on quasi-normal fields is applied. Although the approach can not be carried over straight-forwardly to the case where also the integration on sub-domains is important, that approach has important similarities to the approach presented in this paper.

In the literature of extended finite element methods (XFEM), there exist other approaches which are based on a tesselation approach and aim at improving the resulting piecewise planar approximations by means of a \emph{parametric mapping of the sub-trianguation} \cite{cheng2010higher,dreau2010studied}. These approaches are typically technically involved, especially in more than two dimensions. Moreover, robustness of these approaches is not clear. 

The approach presented in this paper is similar to the approaches in \cite{grande2014highorder,cheng2010higher,dreau2010studied} in that it is also based on a piecewise planar interface which is significantly improved using a parametric mapping. The important difference is, that we consider a \emph{parametric mapping of the underlying mesh} rather than the sub-triangulation or the interface. According to the mapping of the mesh the use of isoparametric finite element is natural. 

At this point, we would also like to mention the paper \cite{Basting2013228} where the construction of a mesh deformation, which is also used in combination with isoparametric finite elements, is specifically designed to align a given mesh to a given interface position. The goal in that paper is to obtain a mesh that is \emph{conforming} to a given interface without changing the mesh topology while keeping the quality of the mesh. Our approach is in a similar virtue. 

\subsection{The concept} \label{sec:intro:mainidea}

In this paper we present a strategy for the numerical approximation of domains implicitly defined by
an approximate signed distance function $\phi$, the level set function.  
The approach is based on the assumption that a numerical strategy to deal with interfaces stemming
from a piecewise linear level set function exactly is available. 

\begin{figure}[h!]
  \vspace*{-0.2cm}
  \begin{center}
    \begin{tikzpicture}[scale=1.0]
      \node[](notcurved)
      {
        \includegraphics[width=0.28\textwidth]{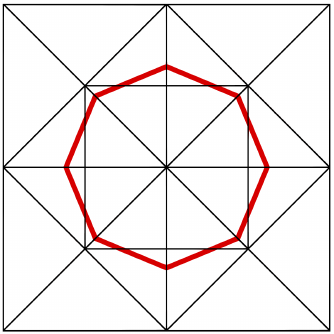} 
      };
      \node[right =2.5cm of notcurved.east, anchor=west](curved)
      {
        \includegraphics[width=0.28\textwidth]{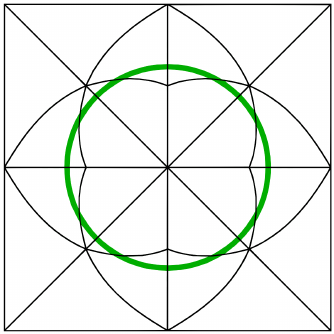} 
      };
      \draw [->] (notcurved.east) to [in=150,out=30] node[above]{$\Psi_h$} (curved.west) ;
    \end{tikzpicture}
  \end{center}
  \vspace*{-0.5cm}
  \caption{Main idea: The piecewise planar interface on the initial mesh is
    mapped close to the desired interface via a mesh transformation $\Psi_h$.}
  \label{fig:idea} 
\end{figure}

In the context of unfitted finite element methods integrals of the form $\int_{S} f \, dx$
have to be computed for $S$ being an interface or a domain which is only implicitly defined through
the level set function $\phi$. Consider the case $S=\Omega_i = \Omega \cap \{ \phi < 0 \}$ for a polygonal
domain $\Omega$ which has a suitable triangulation $\mT$. As an exact evaluation of the integral
$\int_{\Omega_i} f \, dx$ is in general practically not feasible we consider the approximation of
$\Omega_i$ with the domain   
$$
\Omega_{i,h} := \Omega \cap \{ I_h \phi < 0 \}
$$ 
where $I_h$ is the continuous \emph{piecewise linear} interpolation of $\phi$. To $\Omega_{i,h}$
an \emph{explicit representation} can be found on which quadrature rules can be applied element by element:
\begin{equation} \label{eq:pwplan}
\int_{\Omega_i} f \, dx \approx \int_{\Omega_{i,h}} f \, dx \approx \sum_{T\in \mT} \sum_{i}
\omega_i f(x_i)
\end{equation}
Here $\omega_i$ und $x_i$ are the integration weights and points which depend on the cut
configuration on $T$. The accuracy of this approach is
limited by the approximation quality of $\Omega_{i,h}$ which is only second order accurate. 
By a \emph{transformation of the mesh} which is represented by an explicit mapping $\Psi_h$,
parametrized by a finite element function, the piecewise planar interface is mapped approximately
onto the zero level set of a corresponding high order accurate level set function, cf. Figure
\ref{fig:idea}.  

According to this mapping we have that $\Psi_h(\Omega_{i,h})$ is a high order accurate approximation
to $\Omega_i$ which, by construction, still has an explicit representation. The integral can hence be approximated as
\begin{equation} \label{eq:pwplancurved}
\int_{\Omega_i} f \, dx \approx \int_{\Psi_h(\Omega_{i,h})} f \, dx \approx \sum_{T\in \mT} \sum_{i}
\omega_i \, |\mathrm{det}(\Psi_h(x_i))| \, f(\Psi_h(x_i)),
\end{equation}
where the integration weights and points are the same as in \eqref{eq:pwplan}.
In contrast to \eqref{eq:pwplan} the accuracy of the quadrature in \eqref{eq:pwplancurved} is no
longer bounded to second order accuracy but essentially depends on $\Psi_h$. 
The finite element spaces (used in a discretization of an (unfitted) PDE problem) have to be adapted
correspondingly which renders a resulting method an \emph{isoparametric} finite element method: 
Let $W_h$ be the finite element space corresponding to
the piecewise planar interface approximation with $\Gamma_1:=\{I_h \phi = 0\}$. Then, after
transformation with $\Psi_h$, the appropriate isoparametric finite element space is
$$
  \mathcal{W}_h^k := \{ \varphi_h \circ \Psi_h^{-1} | \varphi_h \in W_h \}.
$$
Note that this implies that shape functions are not necessarily polynomials on the mapped domain. 

From a computational point of view only two things change compared to the case of a piecewise planar
interface. First, a suitable (approximated) mapping $\Psi_h$ has to be constructed. A major
contribution of this paper is the discussion of this. Secondly, the deformation has to be considered
in the definition of the shape functions and the quadrature of deformed elements.  
The treatment of the latter aspect is well-known and is not discussed here.

Owing to the combination of a tesselation algorithm for a piecewise planar interface and a
parametric mapping, we obtain the four most important features of this approach: 
\begin{itemize}
\item An explicit high order accurate geometrical approximation of the exact interface. 
\item Guaranteed positiveness of quadrature weights on the interface and in the sub-domains. 
\item Unchanged cut topology compared to the piecewise planar case. 
\item Easy integration into existing codes as the approach builds on a piecewise planar interface.
\end{itemize}
Together, this renders the method highly
accurate, robust, efficient and fairly simple to implement. 
We again note, that the improved
geometry approximation is based on a parametric mapping of the initial mesh, s.t. a combination with
isoparametric finite elements is crucial.

\subsection{Content and structure of the paper}\label{sec:intro:structure}
The main purpose of this study is the presentation of the new approach which allows for an explicite high order geometry approximation of domains which are implicitly defined as level sets. We focus on the discussion of an efficient construction of a suitable mesh deformation and related implementational aspects. Although this discussion goes beyond simple concepts in that it addresses important details of the method, we do not aim at a thorough error analysis, yet. This is topic of ongoing research and will be published in a different paper.

The paper is structured in the following way. The starting point for the method is the ability to
deal with interfaces which are piecewise planar. As we mainly consider simplicial meshes this
piecewise planar interface typically corresponds to a (continuous) piecewise linear level set
function. We consider a corresponding decomposition strategy as given and do not discuss it
here. Instead, we refer to the literature instead, cf. \cite[Chapter 5]{naerland2014geometrychap5}
or \cite[Chapter 4]{lehrenfeld2015diss}.  

The crucial component of the new approach is the construction of an explicit mapping $\Psi_h$ which is suitable for implementation. 
The construction of $\Psi_h$ is discussed in section \ref{sec:discretetrafo} and consists of several steps. 
We first characterize the desired properties of the transformation $\Psi_h$ for a general case. Afterwards we restrict to an important special case with simplicial elements and a given level set function which is an approximate signed distance function. 
For this case we present an explicit construction of a suitable mapping.

Numerical examples demonstrating the quality of the explicit geometrical representation obtained by this new approach are presented in section~\ref{sec:numexgeom} and section~\ref{sec:numexpde}. 
While the examples in section~\ref{sec:numexgeom} focus on the accuracy of the geometry
approximation, in section~\ref{sec:numexpde} a high order unfitted (isoparametric) finite element
formulation for an interface problem is considered which proves the practical use of the method. 

\section{Construction of the mesh transformation $\Psi_h$}  \label{sec:discretetrafo}

We start with the formulation of properties that we demand from a suitable transformation
$\Psi_h$. Let $\Gamma_1$ be a piecewise planar interface and $[V_h^k]^d$ the space of continuous vector-valued piecewise polynomial functions of degree $k$. We seek for a transformation $\Psi_h \in [V_h^k]^d$, s.t. 
\begin{equation} \label{eq:dist}
\Psi_h = \argmin_{\Psi_h \in V^{\ast}} \ \mathrm{dist}(\Psi_h(\Gamma_1),\Gamma) \quad \text{ or at least } \quad \mathrm{dist}(\Psi_h(\Gamma_1),\Gamma) \leq c h^{k+1}
\end{equation}
with $V^\ast \subset [V_h^k]^d$ the subset of admissible transformations which fulfil the following
constraints:
\begin{subequations}
  \begin{itemize}
  \item Homeomorphy: 
    \begin{equation}\label{eq:homeo}
      \Psi_h, \Psi_h^{-1} \in [C^0(\Omega)]^d
    \end{equation}
  \item Shape regularity: Consider an arbitrary element $T$ in the triangulation $\mT$. Let $\Phi_h$
    be the transformation from the reference element $\hat{T}$ to $T=\Phi_h(\hat{T})$. Then the mapped
    element $\mathcal{T} = \Psi_h(T) = (\Psi_h \circ \Phi_h)(\hat{T})$ should also be shape
    regular. Moreover, the shape regularity should be comparable to the shape regularity of $T$. We
    therefore define the transformation $$\Theta_h := \Psi_h \circ \Phi_h : \hat{T} \rightarrow
    \mathcal{T}$$ and demand 
    \begin{equation} \label{eq:cond}
      \kappa(\nabla\Theta_h) \leq C \, \kappa(\nabla \Phi_h).
    \end{equation}
    for a given constant $C>1$.
  \item Locality: 
    \begin{equation}\label{eq:local}
      \Psi_h(x) \neq x \text{ only in the vicinity of the interface } \Gamma_1.
    \end{equation}
  \end{itemize}
\end{subequations}
The last constraint is desirable for efficiency reasons but not crucial. 

Note, that due to the fact that $\Gamma_1$ is already a second order accurate approximation to the
interface, the transformation $\Psi_h$ is essentially only a local high order correction. This is in
contrast to the transformations tailored in \cite{Basting2013228}.
Further note, that $\Psi_h$ is obviously not uniquely defined with the requirements formulated in
\eqref{eq:dist} and \eqref{eq:homeo}-\eqref{eq:local}. Different choices are possible. For the case
of a simplicial mesh we introduce a special choice for $\Psi_h$ in the following.

In order to construct a suitable transformation $\Psi_h$ the above given characterization is too
general to be of direct practical use. 
In the remainder of this study we restrict to simplicial meshes and assume the interface to be
described by an approximate signed distance function $\phi$ which is a continuous piecewise
polynomial function. 
Further, the piecewise planar approximation
$\Gamma_1$ is assumed to be the zero level of the piecewise linear interpolant $I_h \phi$. These
restrictions allow for an \emph{explicit construction} of a function $\Psi_h$ which is of practical
use.  

We derive this in several steps. In section \ref{sec:notation} we introduce notation and assumptions. 
In section \ref{sec:conttrafo} we give an ideal mapping $\Psi$ which fulfils $\Psi(\Gamma_1) =
\Gamma_k:=\{\phi=0\}$ (and hence (\ref{eq:dist})) and is defined pointwise. Together with a suitable
modification (section \ref{sec:psit}) and an approximation of this function in $[V_h^k]^d$ (section
\ref{sec:discretetrafo:fem}) we implement conditions \eqref{eq:dist}, \eqref{eq:homeo} and \eqref{eq:local}. 
In cases where the interface is well resolved by the piecewise linear level set function $I_h \phi$
the transformation $\Psi_h$ is close to the identity and the shape regularity of the transformed
mesh is ensured. Nevertheless, for practical applications a mechanism to ensure robustness even if
the interface is not well resolved becomes necessary.  
The aspect of controlling the shape regularity as in \eqref{eq:cond} is implemented with a
limitation step which is presented and discussed in detail in section
\ref{sec:discretetrafo:shapereg}.  

In section \ref{sec:compasp} we summarize the algorithm to determine $\Psi_h$ and the most important properties of the method from a computational point of view.

\subsection{Notation and assumptions}\label{sec:notation}
We introduce some basic notation and assumptions. 
$\Omega$ is a polygonal domain in $\mathbb{R}^d$, $d=2,3$. 
It is decomposed into a shape regular partition $\mT$ of $\Omega$ consisting of elements $T$ which are \emph{simplices}. 
Inside $\Omega$ an internal interface $\Gamma$ ($\Gamma \cap \partial \Omega = \emptyset$) separates two domains $\Omega_1$ and $\Omega_2$. The interface and the subdomains $\Gamma$, $\Omega_1$ , $\Omega_2$ are implicitely described by a level set function $\phiexact$ with $\Gamma=\{x|\ \phiexact(x) = 0\}$, $\Omega_1=\{x |\ \phiexact(x) < 0\}$ and $\Omega_2=\{x|\ \phiexact(x) > 0\}$. 
The level set function is a signed distance function:
\begin{equation*}
  \phiexact(x) = \left\{ 
    \begin{array}{rc} 
      \text{dist}(x,\Gamma), & x \in \Omega_1, \\
      -\text{dist}(x,\Gamma), & x \in \Omega_2.
    \end{array} 
  \right.
\end{equation*}
We make the following assumptions on the smoothness of the interface $\Gamma$ and the resolution of $\Gamma$ by the triangulation $\mT$. 
\begin{assumption}\label{ass:cmsmooth}
We assume that the interface $\Gamma$ is $C^{m}$-smooth, $m \geq 2$. Then there exists a $\delta > 0$, such that $\phiexact \in H^{m}(\Omega^\delta)$ with
 $\Omega^\delta := \{ x |\  \Vert x -y \Vert_2 \leq \delta \text{ for a } y \in \Gamma\}$.
\end{assumption}
\begin{assumption}\label{ass:p1resolved}
Let $h$ be the maximum mesh size of the triangulation $\mT$. Then we assume that there holds 
$ h < \delta$ with $\delta$ as in assumption \ref{ass:cmsmooth}.
\end{assumption}
Assumption \ref{ass:p1resolved} can always be satisfied for sufficiently small mesh sizes. 
Let $V_h^k$ denote the space of continuous piecewise polynomials of degree $k$. The level set function $\phiexact$ is approximated by a function $\phi \in V_h^k$, $1 \leq k \leq m-1$. We think of $\phi$ as a projection of $\phiexact$ under a suitable projection operator $\Pi^\phi: L^2(\Omega) \rightarrow V_h^k$. In practice one often only knows $\phi$ and not $\phiexact$. We add an assumptions on the approximation error induced by the projection $\Pi^\phi$. 
\begin{assumption}\label{ass:phiexact}
With $\phi = \Pi^\phi \phiexact$ and $\Gamma_k := \{\phi = 0\}$ we assume that there holds \\[-3.5ex]
\begin{subequations}
  \begin{center}
    \begin{tabular}{rr}
      \begin{minipage}{0.5\textwidth}
        \begin{equation}\label{eq:phiapprox1}
          \Vert \phi - \phiexact \Vert_{L^2(\Omega^\delta)} \leq c h^{k+1} \vert \Omega^\delta \vert,
        \end{equation}
      \end{minipage} &
      \begin{minipage}{0.4\textwidth}
        \begin{equation}\label{eq:phiapprox2}
          \Vert \phi - \phiexact \Vert_{\infty, \Omega^\delta} \leq c h^{k+1}, 
        \end{equation}
      \end{minipage}\\
      \begin{minipage}{0.5\textwidth}
        \begin{equation}\label{eq:phiapprox3}
          \Vert \nabla \phi \Vert_{\infty,\Omega^\delta} \leq c \Vert \nabla \phiexact \Vert_{\infty,\Omega^\delta},
        \end{equation}
      \end{minipage}&
      \begin{minipage}{0.4\textwidth}
        \begin{equation}\label{eq:phiapprox4}
          \mathrm{dist}(\Gamma_k,\Gamma) \leq c h^{k+1}.
        \end{equation}
      \end{minipage}
    \end{tabular}
  \end{center}
\end{subequations}
with constants $c$ only depending on $k$ and $\Gamma$.
\end{assumption}
Note that assumption \ref{ass:phiexact} is justified for $\Pi^\phi$ the standard $L^2$ projection on
$\Omega$ or the projection discussed in section \ref{sec:oswald}.

The standard nodal interpolation of $\phi$ into the space of piecewise linears $V_h^1$ is denoted by $I_h \phi$. The zero level of that level set function $\Gamma_1 := \{I_h \phi = 0 \}$ is piecewise planar which allows for an explicit representation of the interface and the sub-domains. $\Gamma_1$ is a second order accurate approximation to $\Gamma$. Let 
$$
\mT^\Gamma := \{ T \in \mT |\, T \cap \Gamma_1 \neq \emptyset \} \text{ and } \Omega^\Gamma := \{ x \in T |\, T \in \mT^\Gamma \}
$$ 
be the set and region of cut elements, respectively. 
Due to assumption \ref{ass:p1resolved} the interface is resolved in the sense that 
$
\Omega^\Gamma \subset \Omega^\delta.
$ 
With the assumed smoothness ($m\geq 2$ in assumption \ref{ass:cmsmooth}) there also holds:
\begin{equation}
\Vert I_h \phi - \phi \Vert_{L^2(\Omega^\Gamma)}
\leq c h^2 \Vert \phiexact \Vert_{H^{2}(\Omega^\Gamma)}. 
\end{equation}
We further introduce the following notation for the extension of cut elements by its neighbors:
$$
\mT^{\Gamma,+} := \{ T \in \mT|\, \partial T \cap \Omega^\Gamma \neq \emptyset \} \text{ and } \Omega^{\Gamma,+} := \{ x \in T |\, T \in  \mT^{\Gamma,+} \}.
$$

\subsection{Ideal transformation $\Psi$} \label{sec:conttrafo}
In this section we define an ideal transformation $\Psi$ which is defined pointwise and is not a
finite element function. We assume the ideal case that $\phi=\phiexact$ and that the assumptions
\ref{ass:cmsmooth} and \ref{ass:p1resolved} are valid. This specifically implies that  $\nabla \phi$
is continuous in the set of cut elements $\Omega^\Gamma$ and $\Gamma = \Gamma_k$. At the end of the
section we discuss the case $\phi \neq \phiexact$.

We ask for $\Psi(\Gamma_1)=\Gamma_k =\Gamma$. This only determines the transformation $\Psi$ 
on $\Gamma_1$. Hence, we use a condition which describes a reasonable extension on $\Omega^\Gamma$. 
Later, in section
\ref{sec:discretetrafo:fem}, the continuous transition to $\Psi_h(x)=x$ away from the
interface is implemented within the approximation in a finite element space.   

\begin{figure}[h!]
  \vspace*{-0.2cm}
  \begin{center}
    \begin{tikzpicture}[scale=1.0]
      \node[](notcurved)
      {
        \includegraphics[width=0.2\textwidth]{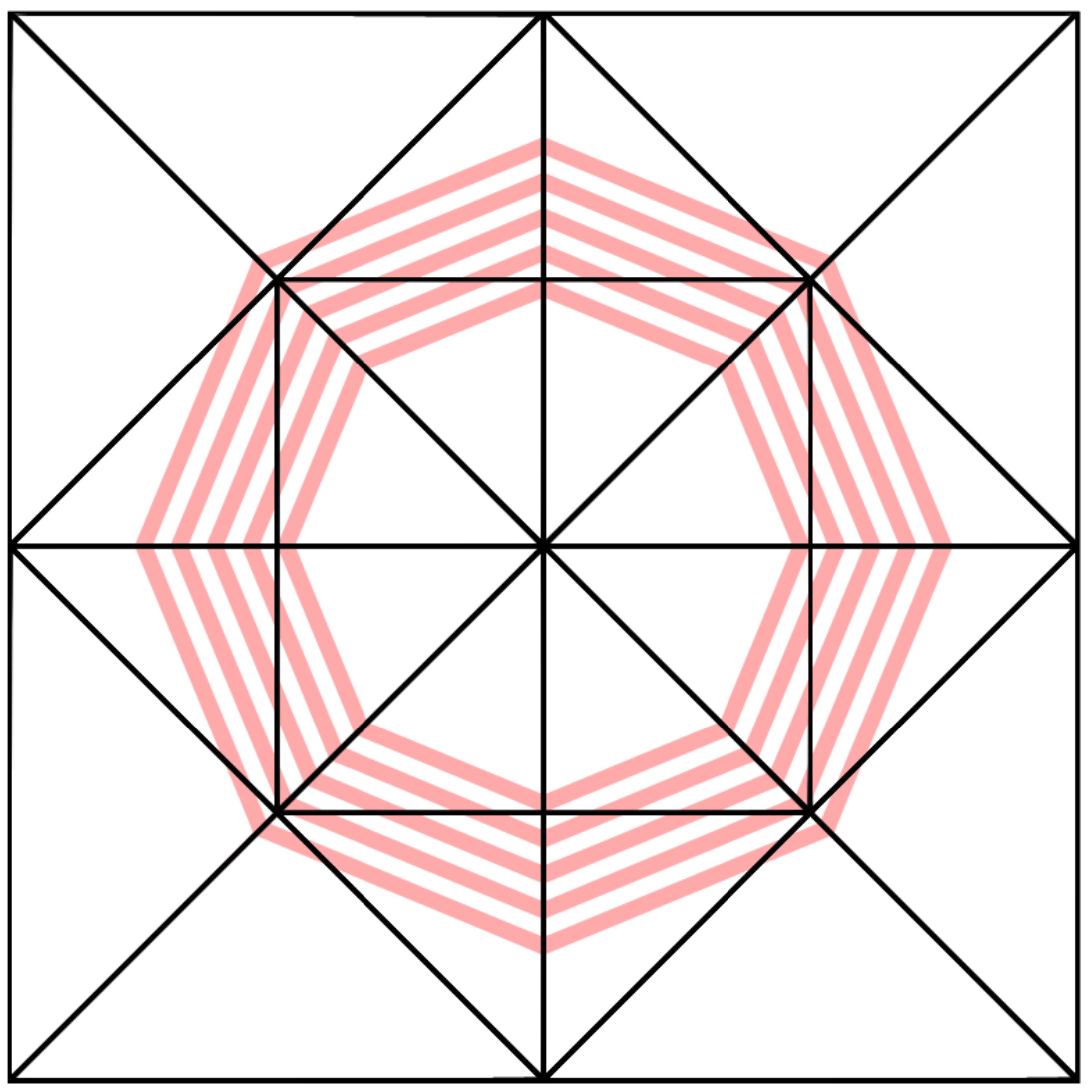} 
      };
      \node[right =2.5cm of notcurved.east, anchor=west](curved)
      {
        \includegraphics[width=0.2\textwidth]{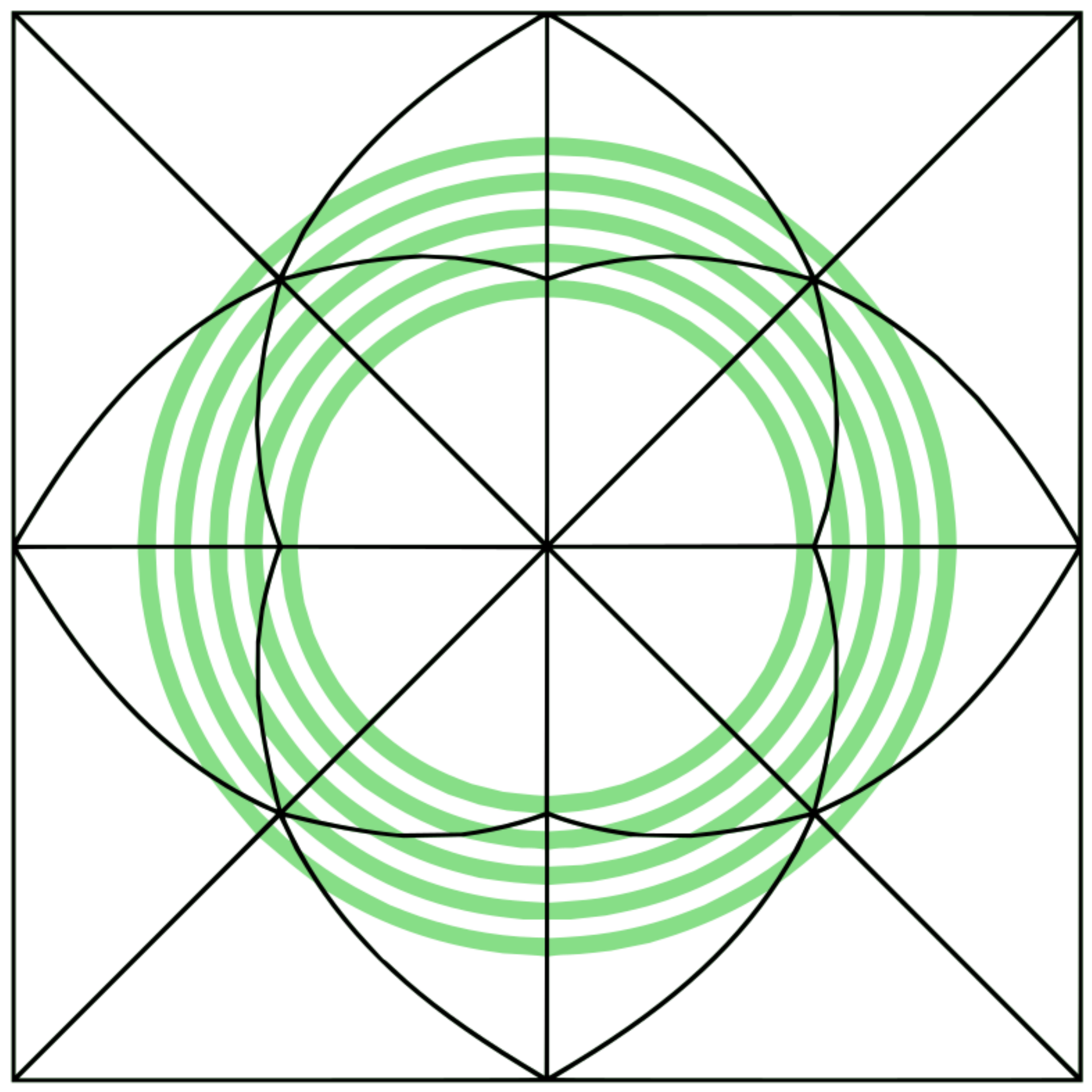} 
      };
      \draw [->] (notcurved.east) to [in=150,out=30] node[above]{$\Psi$} (curved.west) ;
    \end{tikzpicture}
  \end{center}
  \vspace*{-0.5cm}
  \caption{Ideal transformation: Piecewise linear level sets are mapped on exact level sets.}
  \label{fig:idealidea} 
\end{figure}

\paragraph{Definition of the mapping $\Psi$:}
We choose a transformation $\Psi$ which maps iso levels of the piecewise linear approximation of the
level set function $\{ I_h \phi = c \}$ onto the corresponding iso levels of the level set function
$\{ \phi = c \}$, i.e. we ask for
\begin{equation}\label{eq:mapcontours}
 I_h \phi = \phi \circ \Psi \quad \forall \ x \in \Omega^\Gamma.
\end{equation}
Such a transformation is sketched in Figure \ref{fig:idealidea}.
Note that \eqref{eq:mapcontours} implies $\Psi(\Gamma_1) = \Gamma$. The mapping can also be formulated as the following problem: For every point $x \in \Omega^\Gamma$ with (linearly interpolated) level set value $c = I_h \phi(x)$, find a point $x^\ast$ in $\Omega$ such that $\phi(x^\ast) = c$. 

Note that \eqref{eq:mapcontours} does not define $\Psi$ uniquely, yet. In order to determine a unique transformation, for every point $x$ in $\Omega^\Gamma$ we search for the corresponding point $x^\ast \in \{ \phi(x^\ast) = c \}$ with smallest minimal distance to $x$. Due to $x \in \Omega^\delta$ and assumption \ref{ass:cmsmooth} this point is unique.
With the help of the signed distance function $\phi$ we can characterize $x^\ast = x + r s$ using the \emph{search direction} $s = \nabla \phi$ and the distance $r$. 
By specifying the search direction $s(x)$ for each point in $\Omega^\Gamma$, we can define 
$\Psi$ as 
\begin{subequations}
\begin{align}\label{eq:psiastpoint}
  \Psi(x) &:= x + r(x)\, s(x), \ x \in \Omega^\Gamma \text{ with } r(x)\in \rr \text{ such that } \\
\label{eq:linesearch}
  I_h \phi(x) &= \phi(x + r(x)\, s(x)).
\end{align}
\end{subequations}
With $\phi = \phiexact$, we have that $s = \nabla \phiexact$ is continuous in $\Omega^\Gamma$ which results in $\Psi \in [C^0(\Omega^\Gamma)]^d$.
By construction we have $\Psi(x_V)=x_V$ for every vertex $x_V$ in the mesh as there holds $I_h \phi (x_V) = \phi(x_V)$.

\paragraph{The non-ideal case $\phi \neq \phiexact$:}
In practice one mostly has $\phi \neq \phiexact$. Then, typically $s= \nabla \phi$ is not continuous any more and as a result the same holds for a correspondingly defined $\Psi$. There are essentially two ways to deal with this problem. Either continuity of the search directions is restored by means of a projection $\Pi^s$ of $\nabla \phi$ into $[V_h^k]^d$. Possible projections are for instance an $L^2(\Omega^\Gamma)$-projection or the projection discussed in section \ref{sec:oswald}. Alternatively, one could stick with the discontinous search direction $s = \nabla \phi$ accepting the resulting discontinuous $\Psi$. This is acceptable as finally, for a realization of the transformation, $\Psi$ is projected into a continuous finite element space with a projector $\Pi^\Psi$ similar or identical to $\Pi^s$. 
In our experience both approaches give similar results. Nevertheless, in the remainder of the paper we consider the use $s = \Pi^s \nabla \phi$ as it simplifies the discussion of the method. 
We also mention that it is not necessary to have a high order accurate approximation of $\nabla
\phi$ for the search direction. In the numerical examples in section \ref{sec:numexgeom} and section
\ref{sec:numexpde} a \emph{continuous piecewise linear} search direction (using the projection
discussed in section \ref{sec:oswald}) is used.

\begin{figure}[h!]
  \vspace*{-0.2cm}
  \begin{center}
    \includegraphics[width=0.5\textwidth]{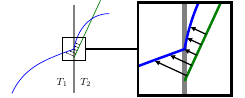} 
  \end{center}
  \vspace*{-0.5cm}
  \caption{Sketch of a situation where the ideal transformation $\Psi$ is not smooth within a given
    element ($T_2$). Points on the piecewise planar level set (green) are mapped on corresponding points on
    the corresponding level set of $\phi$ (blue). $\Psi$ inherits the missing regularity of $\phi$
    with a shift in locations: kinks across element interfaces in $\phi$ result in kinks of $\Psi$
    within elements. }
  \label{fig:nonsmooth} 
\end{figure}

Another important aspect of the non-ideal case $\phi \neq \phiexact$ is the fact, that the
transformation $\Psi$ is typically not smooth (only Lipschitz continuous) within each element $T$ as
$\phi$ typically has discontinuities in all derivatives across element interfaces. A typical situation
leading to a non-smooth $\Psi$ is sketched in Figure \ref{fig:nonsmooth}. In section \ref{sec:psit}
we propose a modification which (alongside its other benefits) results in a transformation which is smooth
within all elements.

\paragraph{Evaluation of $\Psi$:} \label{sec:discretetrafo:newton}
At each point where the transformation $\Psi$ has to be evaluated, the solution of problem \eqref{eq:linesearch} for a fixed $x$ has to be implemented. Let $T$ be the simplex, such that $x \in T$.
The search direction $s(x)$ and the value $(I_h \phi)(x)$ can easily be evaluated at the given point. This information is available element-local. In order to find the corresponding projected point $x^\ast$, it seems natural to try the following Newton iteration:
\begin{equation}\label{eq:newton}
 x^0 = x, \quad x^{k+1} = x^k - \frac{(I_h \phi)(x) - \phi(x^k)}{\nabla \phi(x^k) \cdot s(x)} \cdot s(x)
\end{equation}
This procedure has two disadvantages. First, $\phi$ may have kinks at element interfaces (as it is
only piecewise polynomial) which may complicate the solution of the non-linear problem. This problem
can be tackled by replacing the Newton-iteration by more sophisticated methods. More importantly,
for $x^0 \in T$ it can happen that $x^k \notin T$, such that evaluations of $\phi$ need to be
carried out in a neighborhood of element $T$. The evaluation of finite element functions on neighbor elements is often troublesome, esspecially if one is concerned with an efficient parallelization. We propose a modification which solves both problems in the next section. 


\subsection{Localized (piecewise smooth) transformation $\Psi_\mT$}\label{sec:psit}
We modify the previously discussed construction of the deformation in order to render an
implementation less complicated and more efficient.
 We achieve this with a localization of the computations on elements. That means that we avoid
 evaluating finite element functions from neighbor elements. To this end we first introduce a
 discontinuous (across element interfaces) transformation $\Psi_\mT$ in section
 \ref{sec:psit} which is afterwards projected into the space of continuous finite element functions
 with a projection such as the one discussed in section \ref{sec:oswald}.   

\begin{figure}[h!]
  \vspace*{-0.2cm}
  \begin{center}
    \includegraphics[width=0.5\textwidth]{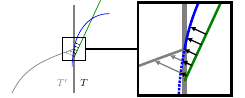} 
  \end{center}
  \vspace*{-0.5cm}
  \caption{Sketch of the same configuration as in Figure \ref{fig:nonsmooth}, this time for the localized transformation
    $\Psi_T$. Points on the piecewise planar level set (green) are mapped on corresponding points on
    the corresponding \emph{locally extended} (dotted) level set of $\phi$ (blue) resulting in a smooth
    function $\Psi_T$ within each element $T$.}
  \label{fig:nonsmoothfixed} 
\end{figure}

We define a modified transformation $\Psi_\mT$ by its element-wise contributions $\Psi_T$ and only
consider elements $T \in \mT^\Gamma$. On each element $T$ the level set $\phi$ is a polynomial of order $k$. Let $\mathcal{E}_T$ be the extension operator so that $\mathcal{E}_T(\phi)$ denotes the polynomial on $\rr^d$ the restriction of which on $T$ coincides again with $\phi$. As $\phi$ is an accurate approximation of a smooth function, $\mathcal{E}_T(\phi)$ is a good approximation to $\phi$ also in the neighborhood of $T$. 
Replacing $\phi$ with $\mathcal{E}_T(\phi)$ in \eqref{eq:linesearch}, results in the definition of $\Psi_T$:
\begin{subequations}
\begin{align}\label{eq:psimod}
  \Psi_T(x) &:= x + r_T(x)\, s(x), \ x \in T \text{ with } r_T(x)\in \rr \text{ such that } \\
\label{eq:linesearchmod}
  I_h \phi(x) &= (\mathcal{E}_T (\phi))(x + r_T(x)\, s(x)).
\end{align}
\end{subequations}
As the function to evaluate in \eqref{eq:linesearchmod} is now a polynomial on $\rr^d$ the Newton iteration
\begin{equation}\label{eq:modnewton}
  x^{k+1} = x^k - \frac{(I_h \phi)(x) - (\mathcal{E}_T(\phi))(x^k)}{\nabla (\mathcal{E}_T(\phi))(x^k) \cdot s(x)} \cdot s(x)
\end{equation}
is much simpler to evaluate and should converge even faster. For moderate order $k$ even exact root finding algorithms could be applied.
Another advantage of this approach is that we can guarantee that the resulting transformation is
smooth within each element, cf. Figure \ref{fig:nonsmoothfixed}. 

We briefly explain why we expect this to be a valid approximation of the original problem. To this end, we assume that the interface is resolved. In that case $\mathcal{E}_T(\phi)$ only has to be evaluated in a small neighborhood of $T$, denoted as $\omega(T)$. In that neighborhood a simple triangle inequality gives
\begin{equation}\label{eq:estextens}
 \Vert \mathcal{E}_T(\phi) - \phi \Vert_{\infty,\omega(T)} \leq  \Vert \mathcal{E}_T(\phi) - \phiexact \Vert_{\infty,\omega(T)} +  \Vert \phi - \phiexact \Vert_{\infty,\omega(T)} 
\end{equation}
where the latter term can be bounded in $\mathcal{O}(h^{k+1})$ due to assumption \ref{ass:phiexact} and the first term on the right hand side can also be bounded in $\mathcal{O}(h^{k+1})$ with a Bramble-Hilbert argument.

$\Psi_\mT$ is discontinuous across element boundaries. However, the jumps over element boundaries should be in the order of $\mathcal{O}(h^{k+1})$. Again by a projection $\Pi^\Psi$ we make sure that the final discrete transformation $\Psi_h$ is a continuous finite element function. A simple and efficient choice for $\Pi^\Psi$ is discussed in section \ref{sec:oswald}.

\subsection{Approximation with finite elements} \label{sec:discretetrafo:fem}
Until now we have discussed how to define and evaluate $\Psi$ in $\Omega^\Gamma$. Our goal, however,
is the construction of a finite element function $\Psi_h \in [V_h^k]^d$ which fulfils condition
\eqref{eq:dist}, but furthermore guarantees \eqref{eq:homeo} and \eqref{eq:local}.

We write $\Psi_h \in [V_h^k]^d$ as $\Psi_h = x + d_h$ with $d_h \in [V_h^k]^d$ the function describing the deformation of the mesh. Now we can split
$$
d_h = d_h^\Gamma + d_h^\Omega = \sum_{i=1}^{N_\Gamma} u_i \varphi^\Gamma_i(x) + \sum_{j=1}^{N_\Omega} v_j \varphi^\Omega_i(x)
$$
with $u_i \in \rr^d$ the degrees of freedom corresponding to the basis functions $\varphi_i^\Gamma$ which are located at the interface, $\mathrm{supp}(\varphi^\Gamma_i) \cap \Omega^\Gamma \neq \emptyset$, and $v_j$ the remaining degrees of freedom ($\mathrm{supp}(\varphi^\Omega_j) \cap \Omega^\Gamma = \emptyset$). By construction we then have $\Psi_h|_{\Omega^\Gamma} = (x + d_h^\Gamma)|_{\Omega^\Gamma}$ and $d_h^\Omega|_{\Omega^\Gamma} = 0$, i.e. $d_h^\Omega$ has no contribution in $\Omega^\Gamma$.

This splitting decomposes the problem again into two parts. First, $d_h^\Gamma$ has to be chosen
such that $\Psi_h$ is a good approximation to $\Psi$ in $\Omega^\Gamma$. Secondly, $d_h^\Omega$ has
to be chosen such that $\Psi_h(x) = x$ in $\Omega \setminus \Omega^{\Gamma,+}$. The second problem
immediately suggests to choose $d_h^\Omega = 0$. Hence, \eqref{eq:homeo} and \eqref{eq:local} are
fulfilled and we only have to consider the approximation in \eqref{eq:dist}. On that account we
apply a projection $\Pi^\Psi$ of $\Psi_{\mT}$ (or $\Psi$) into $[V_h^k]^d$. This projection could be
a nodal interpolation (only if $\Psi$ is used), an $L^2(\Omega^\Gamma)$-projection or the projection
in section \ref{sec:oswald} and we expect the following assumption to be valid. 
\begin{assumption}\label{ass:psiapprox}
Consider $\Psi^\ast \in \{ \Psi, \Psi_{\mT} \}$.  With $\Psi_h = \Pi^\Psi \Psi^\ast$ there holds 
\begin{equation}
  \Vert \Psi_h - \Psi^\ast \Vert_{L^2(\Omega^\Gamma)} \leq c h^{k+1} \Vert \Psi^\ast \Vert_{H^{k+1}(\mT^\Gamma)} \leq c h^{k+1} \vert \Omega^\Gamma \vert
\end{equation}
where the last inequality holds due to $\Psi^\ast(x) \stackrel{h\rightarrow 0}{\longrightarrow} x$. Further there holds
\begin{equation}\label{eq:psiinfty}
  \Vert \Psi_h - \Psi^\ast \Vert_{\infty,\Omega^\Gamma} \leq c h^{k+1}.
\end{equation}
\end{assumption}
A rigorous justification of this assumption (for both, $\Psi$ and $\Psi_{\mT}$) requires further
analysis which is topic of ongoing research. 

Consider $\Psi_h = \Pi^\Psi \Psi$. Defining $\Gamma_h := \Psi_h ( \Gamma_1 )$ we expect the
following error estimate to hold:
\begin{align}
\mathrm{dist}(\Gamma_h, \Gamma) & \leq \Vert \phiexact \Vert_{\infty,\Gamma_h} = \Vert \phiexact \circ \Psi_h \Vert_{\infty,\Gamma_1} = \Vert \phiexact \circ \Psi_h - I_h \phi \Vert_{\infty,\Gamma_1} \nonumber \\[0.75ex]
 & \leq \Vert \phiexact \! \circ \! \Psi_h\! -\! \phi \circ \! \Psi_h \Vert_{\infty,\Gamma_1} \! + \Vert \phi \circ \!\Psi\! -\! \phi \circ \!\Psi_h \Vert_{\infty,\Gamma_1} \!+\! \underbrace{\Vert \phi \circ \!\Psi - I_h \phi \Vert_{\infty,\Gamma_1}}_{=0} \label{eq:errest}\\[-2.5ex]
& \leq \underbrace{\Vert \phiexact - \phi \Vert_{\infty,\Omega^\Gamma}}_{\leq ch^{k+1}}
 + \underbrace{\Vert \nabla \phi \Vert_{\infty, \Omega^\Gamma}}_{\leq c \Vert \nabla \phiexact \Vert_{\infty,\Omega^\Gamma}}
\underbrace{\Vert \Psi - \Psi_h \Vert_{\infty, \Gamma_1}}_{\leq ch^{k+1}} \leq c h^{k+1} \nonumber
\end{align}
where in the last step we used assumptions \ref{ass:phiexact} and \ref{ass:psiapprox}.
For the case $\Psi_h = \Pi^\Psi \Psi_{\mT}$ the estimate still applies (replace $\Psi$ with
$\Psi_{\mT}$) except for the term $\Vert \phi \circ \!\Psi_{\mT} - I_h \phi \Vert_{\infty,\Gamma_1}$
in \eqref{eq:errest} which is not exactly zero but can be estimated as follows:

\begin{align*}
  \Vert \phi \circ \!\Psi_{\mT} - I_h \phi \Vert_{\infty,\Gamma_1} 
  & \leq 
    \max_{T \in \mT^\Gamma} \big\{ \underbrace{\Vert \mathcal{E}_T(\phi) \circ \!\Psi_{\mT} - I_h \phi
    \Vert_{\infty,\Gamma_1\cap T} }_{=0} + \Vert \mathcal{E}_T(\phi) \circ \!\Psi_{\mT} - \phi \circ \!\Psi_{\mT}
    \Vert_{\infty,T}  \big\}\\ 
  & \leq
    \max_{T \in \mT^\Gamma} \Vert \mathcal{E}_T(\phi) - \phi
    \Vert_{\infty,\omega(T)} \stackrel{\eqref{eq:estextens}}{\leq} c h^{k+1}\\ 
\end{align*}

\subsection{An Oswald-type projection}\label{sec:oswald}
At several occasions throughout this paper we mentioned a generic projection operator $\Pi: [L^2(\Omega^\ast)]^n \rightarrow [V_h^k]^n$ with $\Omega^\ast = \Omega$ or $\Omega^\ast = \Omega^\Gamma$ and $n \in \{1,d\}$, for instance $\Pi^\phi$, $\Pi^s$ and $\Pi^\Psi$. An obvious choice for all these projection is a simple $L^2(\Omega^\ast)$ projection. However, to avoid the solution of global linear systems it is worth mentioning a simple and more efficient alternative which we also used in the numerical examples.
For ease of presentation we only consider the case where we seek a projection $\Pi:
L^2(\Omega^\Gamma) \rightarrow V_h^{k}$, i.e. the scalar case on the interface region. The
translation to $\Omega$ and/or vector-valued functions is then obvious. 

The Oswald-type interpolation operator consists of two steps: A projection into a discontinuous finite element space and an averaging into $V_h^k$. 
Let $V_h^{k,\text{disc}} := \{u|_T \in \mathcal{P}^k(T) | \, T \in \mT^\Gamma\}$ be the space of piecewise polynomials which are discontinuous across element boundaries. Due to the missing continuity-restriction, the standard $L^2$-projection $\Pi^\text{disc}: L^2(\Omega) \rightarrow V_h^{k,\text{disc}}$ can be computed in an element-by-element fashion and can be implemented very efficiently. Next, we define a simple averaging operator $\Pi^{\text{av}}: V_h^{k,\text{disc}} \rightarrow V_h^{k}$ which is a high order version of the \emph{Oswald interpolation operator} \cite{oswald}. Let $\omega_i := \{ T \in \mT^\Gamma |\, \mathrm{supp}(\varphi_i) \cap T \neq \emptyset \}$ be the set of elements where the $i$th basis function $\varphi_i$ of $V_h^k$ is supported in $\Omega^\Gamma$. Further, let $u_{i,T}$ be the coefficients such that 
$$
u = \sum_{i} \sum_{T \in \omega_i} u_{i,T} \, \varphi_i|_T, \quad u \in V_h^{k,\text{disc}}.
$$
We then define
$
  \Pi^{\text{av}}(u) := \sum_i w_i \varphi_i
$
\begin{equation*}
\text{ with } w_i = \frac{1}{\# \omega_i}  \sum_{T\in\omega_i} u_{i,T} \text{ if } \mathrm{supp}(\varphi_i) \cap \Omega^\Gamma \neq \emptyset \text{ and } w_i = 0 \text{ otherwise.}
\end{equation*}
Together, we obtain the projection operator $\Pi^o: L^2(\Omega) \rightarrow V_h^k$ with $\Pi^o = \Pi^{\text{av}} \circ \Pi^{\text{disc}}$. We note that in an implementation of $\Pi^o$ the finite element space $V_h^{k,\text{disc}}$ does not need to be constructed explicitly.

\subsection{Shape regularity} \label{sec:discretetrafo:shapereg}
So far we assumed that the interface is smooth and well resolved. An approach which aims to be of
practical use should however be capable of dealing with non-smooth or underresolved interfaces. In
this section we characterize a simple sufficient condition to ensure shape regularity. Based on this
we discuss two different cases. First, we assume the nice case, where the interface is well-resolved
and smooth and show that shape regularity is not an issue. Secondly, we discuss the case where the
interface is not resolved and propose a simple strategy to ensure shape regularity. 

Before we discuss the details of both cases, we introduce transformations between the reference
simplex, the simplex before and after the transformation with $\Psi_h$ and an intermediate curved
simplex.  

\paragraph{Simplex transformations:}  
A planar simplex can be characterized as a reference simplex under an affine transformation. The mesh obtained after the transformation with $\Psi_h$ can be interpreted as a collection of simplices under the concatenation of two transformations: an affine one and another one parametrized by a finite element function. This characterization facilitates the investigation of the shape regularity of the transformed elements. We therefore introduce notation for the (curved and planar) simplices and the transformations relating them. In figure \ref{fig:trafos} the simplices and transformations are also depicted. 

By $\hat{T}$ and $\hat{x}$ we denote the reference simplex and its coordinates. The affine mapping ${\Phi}_h$ transforms the reference simplex to the planar simplex ${T} = {\Phi}_h(\hat{T})$. 
\begin{subequations}
\begin{equation}
{x} = {\Phi}(\hat{x}) = A \cdot \hat{x} + x_0
\end{equation}
The final (curved) element $\mathcal{T}$ is obtained after the finite element transformation of
${T}$ by $\Psi_h$. The shape functions (corresponding to $\Psi_h$) are typically defined w.r.t. to the reference simplex $\hat{T}$, s.t. we have 
\begin{equation}
y = \Psi_h ({x}) = {x} + \sum_{i} d_i\, (\varphi_i \circ {\Phi}_h^{-1}) ({x}) \overset{x=\Phi_h(\hat{x})}{=} x_0 + A \cdot \hat{x} + \sum_i
d_i \cdot \varphi_i(\hat{x}).
\end{equation}
\begin{figure}[] 
\vspace*{-0.2cm}
  \begin{center}\includegraphics[width=0.8\textwidth]{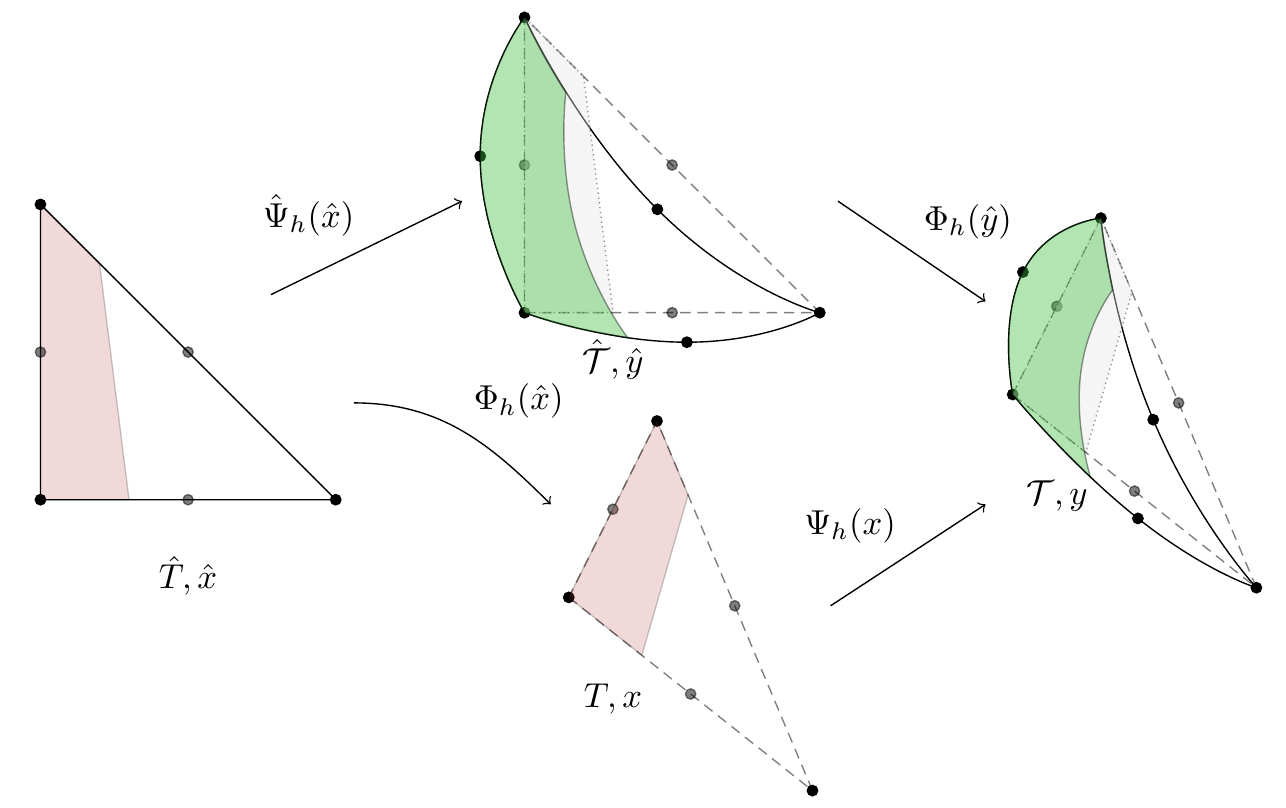}\end{center}
\vspace*{-0.6cm}
  \caption{Reference and transformed domains before and after curving and before and after affine transformation to the physical domain.}
\label{fig:trafos}
\end{figure}
Thus, we have $\mathcal{T} = \Theta_h(\hat{T})$ with $\Theta_h = \Psi_h \circ \Phi_h$. We can also characterize $\Theta_h$ as a mapping which first applies the non-planar deformation and the affine mapping afterwards. We therefore introduce the mapping $\hat{\Psi}_h$, s.t.
\begin{equation}
\Theta_h = (\Psi_h \circ \Phi_h) = (\Phi_h \circ \hat{\Psi}_h)
\end{equation}
with 
\begin{equation}
\hat{y} = \hat{\Psi}_h(\hat{x}) = \hat{x} + \sum_i \hat{d}_i \cdot \varphi_i(\hat{x}), \quad \hat{d}_i = A^{-1} \cdot {d}_i.
\end{equation}
\end{subequations}
Let $\hat{x}_V$ and $x_V$ be the coordinates of the vertices of $\hat{T}$ and $T$. Then $\hat{\Psi}_h(\hat{x}_V) = \hat{x}_V$ and $\Psi_h(x_V) = x_V$.

From the previous observations and we can conclude that
\begin{equation}
  \nabla \Theta = \nabla \Phi_h + \sum_i d_i \nabla \varphi_i = \nabla \Phi_h \cdot \nabla \hat{\Psi}_h, \quad \nabla \Phi_h = A.
\end{equation}
Hence, if the transformation $\hat{\Psi}_h: \hat{T} \rightarrow \hat{\mathcal{T}}$ (curving of the reference element) and the transformation ${\Phi}_h : {T} \rightarrow T$ (affine transformation) are well-behaved, the same also follows for $\Theta: \hat{T} \rightarrow \mathcal{T}$.

As a measure to quantify the shape regularity we consider the relative spectral condition number of
the Jacobian. We have the simple estimate
\begin{equation}
\kappa ( \nabla \Theta ) \leq \kappa(\nabla \Phi_h) \cdot \kappa(\nabla \hat{\Psi}_h).
\end{equation}
We may assume that the initial mesh is shape regular such that $\kappa(\nabla \Phi_h) \leq C$ for
reasonable moderate number $C$. In order to fulfil \eqref{eq:cond} it remains to ask for the impact of the non-linear transformation, i.e. $\kappa(\nabla {\hat{\Psi}_h})$.

We have $\nabla \hat{\Psi}_h = I + \nabla \hat{d}_h$ with $\hat{d}_h(\hat{x})$ the deformation of $\hat{\mathcal{T}}$ w.r.t. the reference simplex $\hat{T}$. As $\hat{d}_h$ is a polynomial on the reference triangle, we can bound its gradient
\begin{equation} \label{eq:inverse}
  \Vert \nabla \hat{d}_h \Vert_{\infty,\hat{T}} \leq c_k \Vert \hat{d}_h \Vert_{\infty,\hat{T}}
\end{equation}
with $\Vert q \Vert_{\infty,\hat{T}} = \max_{x\in \hat{T}} \Vert q (x) \Vert_2$ for a vector-valued function $q$
and $\Vert Q \Vert_{\infty,\hat{T}} = \max_{x\in \hat{T}} \Vert Q (x) \Vert_F$ for a matrix-valued function $Q$. The constant $c_k$ only depends on the polynomial degree. Due to the choice of norms, it is sufficient to have \eqref{eq:inverse} for a scalar-valued polynomial.

A simple sufficient condition for the shape regularity, formulated as  $\kappa(\nabla \hat{\Psi}_h)
\leq C$ for a chosen constant $C > 1$, can be derived as follows. With elementary
calculations we have
\begin{align*}
 & \Vert \nabla \hat{\Psi}_h \Vert_2 && \leq 1 + \Vert \nabla \hat{d}_h \Vert_2 && \leq  1 + \Vert \nabla
  \hat{d}_h \Vert_F && \leq  1 + \Vert \nabla \hat{d}_h \Vert_{\infty,\hat{T}} \\ 
 & \Vert (\nabla \hat{\Psi}_h)^{-1} \Vert_2 && \leq (1 - \Vert \nabla \hat{d}_h \Vert_2)^{-1} && \leq (1 - \Vert \nabla \hat{d}_h \Vert_F)^{-1} && \leq 
 (1 - \Vert \nabla \hat{d}_h \Vert_{\infty,\hat{T}})^{-1}.
\end{align*}
With \eqref{eq:inverse} we can easily conclude that 
\begin{equation}\label{eq:suff}
 \Vert \hat{d}_h \Vert_{\infty,\hat{T}} \leq \frac{C - 1 }{c_k (C + 1) } 
\end{equation}
is a sufficient condition for
$$
 \kappa( \nabla \hat{\Psi}_h) \leq C.
$$
\paragraph{The resolved case:}
Due to $ \Vert I_h \phi - \phi \Vert_{L^2(\Omega^\Gamma)} \leq c h^2$ there holds $\mathrm{dist}(\Gamma_1,\Gamma_k) \leq ch^2$ and we can conclude that $\Vert \hat{d}_h \Vert_{\infty,\hat{T}} \leq c h$. This implies that the deformation gets arbitrary small which implies \eqref{eq:suff} for sufficiently small $h$. We can thus conclude, that the shape regularity is easily provided for sufficiently small mesh sizes, i.e. for a sufficiently fine resolution of the interface by the mesh. 

\paragraph{The under-resolved case:}
The considered construction of the transformation may lead to undesirable deformation, e.g. self-intersecting geometries and arbitrary small angles if the interface is not sufficiently well resolve. This may lead to a blow up in $\kappa(\nabla \hat{\Psi}_h)$ or even a singular $\nabla \hat{\Psi}_h$ at some points. 

To avoid these situations, we introduce a ``barrier step'' which garantuees the quality of the resulting mesh. The idea is simple. 
We replace the deformation $d=\Psi-x$ (or $d=\Psi_{\mT} - x$) with a \emph{limited deformation}
\begin{equation}\label{eq:barrier}
\bar{d} = \min\{ \gamma h \cdot {s}/{\Vert s \Vert}, d \}.
\end{equation}
On a fixed element $T$, we then have $\hat{d}_h = (\nabla \Phi_h)^{-1} \Pi^\Psi \bar{d}$ and hence
\begin{equation}
  \Vert \hat{d}_h \Vert_{\infty} \leq c_T c_\pi \gamma
\end{equation}
with $c_\pi$ the $L^\infty$ stability constant of the projection $\Pi^\Psi$ and $c_T = h \Vert (\nabla \Phi_h)^{-1} \Vert_{\infty,T}$ a constant depending on the shape regularity of $T$ and quasi-uniformity of $\mT^\Gamma$.

For sufficiently small parameter $\gamma$ we can then ensure condition \eqref{eq:suff} and get shape regularity independent of the resolution of the interface. The price, however, is a reduction to an essentially only second order accurate approximation if the limitation step has to be used ($\bar{d} \neq d$). 

In cases where the interface is resolved the limitation is not needed
which leads to $ \overline{d}=d$. In the worst case, where the interface is not resolve, one
essentially recovers the quality of the piecewise linear approximation.  
In both cases the deformed elements are shape regular which results in a robust method.

\subsection{Summary of computational aspects} \label{sec:compasp}
To clarify on the (low) computational complexity of the resulting scheme we briefly summarize the algorithmic structure for the computation of $\Psi_h$. We restrict to the modified version of the algorithm, cf. section \ref{sec:psit}, which is also used in the numerical examples. 
To this end we determine the coefficients $d_i \in \rr^d$ of the representation $\Psi_h(x) = x + \sum_i d_i \varphi_i(x),~d_i \in \rr^d$ where $\varphi_i$ are the basis functions of $V_h^k$:
\begin{enumerate}
\item Set coefficients and counter (for later averaging) to zero: \\
  $d_i=0$ and $n_i=0$, $i=1,..,\mathrm{dim}(V_h^k)$.
\item Compute element contributions ($\Pi^{\text{disc}} \Psi_{\mT}$):\\
  Loop over elements $T \in \mT^\Gamma$:
  \begin{enumerate}
  \item Determine the $\Psi_{h}|_T \in [\mathcal{P}^k(T)]^d$ by solution of the local $L^2(T)$ problem:
    $$
    \int_{T} \Psi_{h}|_T \cdot v \, dx = \int_T \Psi_T \cdot v \, dx \quad \forall \, v \in [\mathcal{P}^k(T)]^d
    $$
    For each integration point $x_j$ this involves:
    \begin{itemize}
    \item The evaluation of $\Psi_T(x_j)$:\\
      This requires the solution of \eqref{eq:linesearchmod} by the Newton iteration \eqref{eq:modnewton}.
    \item The limitation step:\\
      If the deformation $\Psi_T(x_j)\! -\! x_j$ is too large, it is restricted according to \eqref{eq:barrier}.
    \end{itemize}
  \item Add element contribution to global coefficient vectors (for all d.o.f. $i$ of $T$):
    \begin{itemize}
    \item $d_i \leftarrow d_i + d_i^T$ where $d_i^T$ are the basis coefficients corresponding to $\Psi_{h}|_T$
    \item $n_i \leftarrow n_i + 1$
    \end{itemize}
  \end{enumerate}
\item Average the element contributions ($\Pi^{\text{av}}$): \\
  $d_i \leftarrow d_i/n_i$ if $n_i > 0$
\end{enumerate}
In this overview we neglected the computation of the search direction $s=\Pi^s \nabla \phi$. The corresponding algorithm goes the same lines with the only difference that the evalution of the right hand side functional in the (element-)local $L^2$ projection is much simpler.

The resulting deformation $\Psi_h$ has the following properties:
\begin{itemize}
\item Only elements in the neighborhood of $\Gamma$ are deformed. In $\Omega \setminus \Omega^{\Gamma,+}$ the mesh stays unchanged. 
\item Where the resolution of the interface is sufficiently high, the interface can be resolved with high order accuracy.
\item In any case shape regularity can be guaranteed. 
\end{itemize}

\section{Numerical examples I: Geometry approximation} \label{sec:numexgeom}

\subsection*{Preliminaries}
In the following we consider complex geometries described by level set functions and their approximation with the approach proposed before. For the barrier parameter we choose $\gamma = 0.1$ (independent of the degree $k$). To solve the one-dimensional nonlinear problem in \eqref{eq:linesearch} we choose a (relative) tolerance of $10^{-14}$ which resulted in only of $2-3$ iterations for each point. 

For the quasi-normal field $s$ we use a \emph{continuous piecewise linear} field. In all examples we use the transformation $\Psi_h = \Pi^o (\Psi_\mT)$, cf. section \ref{sec:psit}.
The implementation is included in an add-on library \texttt{ngsxfem} to the finite element library \texttt{NGSolve} \cite{schoeberl2014cpp11}.

We are interested in the geometrical error between the discrete interface $\Gamma_h = \Psi_h(\Gamma_1)$ (with an explicit description) and the ideal interface $\Gamma = \{ \phiexact(x) = 0 \}$ (with an implicit description). Note that the error consists of two approximations: The approximation of the level set function $\phi \approx \phiexact$ and the approximation of $\Gamma_k = \{\phi = 0\}$ with $\Psi_h(\Gamma_1)$. The error we are interested in is the maximum distance between the discrete and the exact interface, $\text{dist}({\Gamma}_h,\Gamma)$. 

\subsection{A two-dimensional example: The flower-shape}\label{sec:flower}

\begin{figure}[h!]
  \begin{center}
  \includegraphics[trim=7cm 0cm 8cm 0cm, clip=true, width=0.38\textwidth]{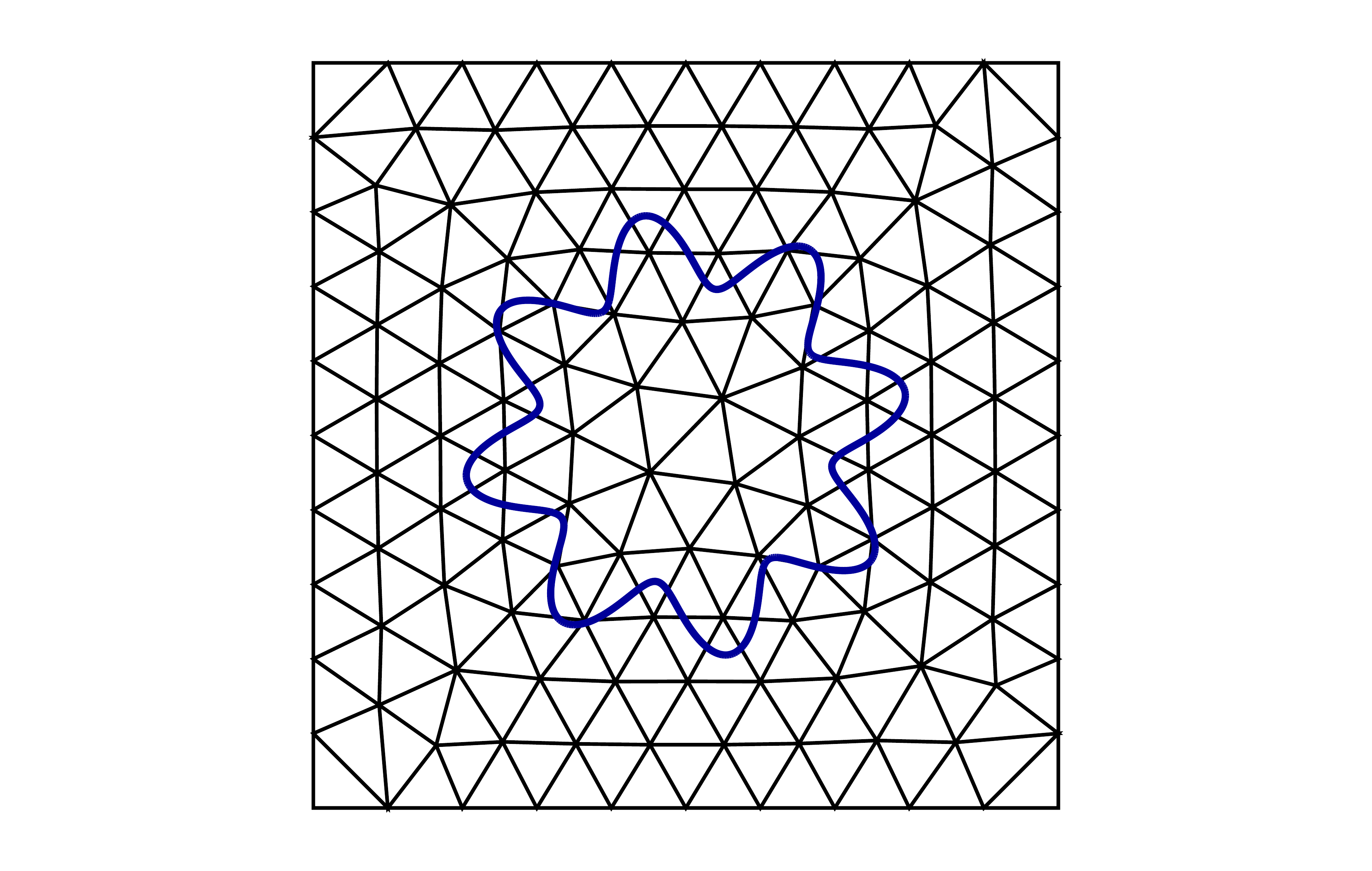}
\includegraphics[trim=0.25cm 0cm 0.25cm 0cm, clip=true, width=0.60\textwidth]{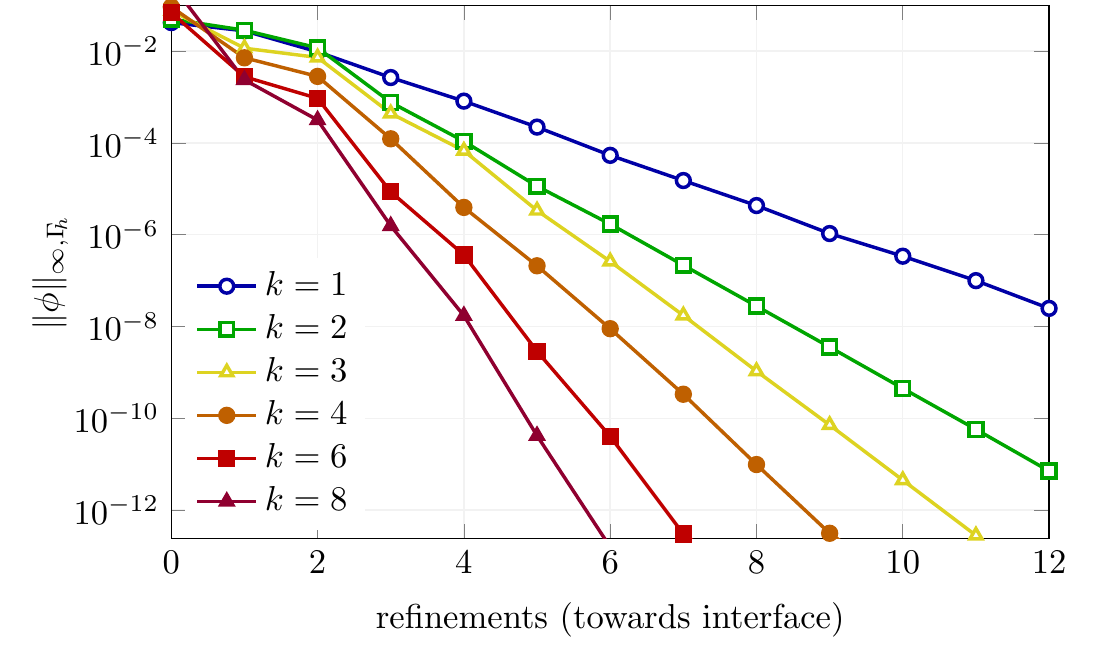} 
\end{center}
  \vspace*{-0.75cm}
  \caption{Mesh and convergence behavior of geometrical error for the flower example.}
\label{fig:flower}
\end{figure}

The geometry in this example is inspired by the one in \cite{li95}. Due to its shape we call it the ``flower-shape'' example. The domain is a square domain $\Omega=[-1,1]^2$ and the interface is given as 
\begin{equation}
\Gamma := \{ (x_1,x_2) |\ x_1=R(\theta) \cos(\theta), x_2=R(\theta) \sin(\theta), \theta \in [0,2\pi]\}
\end{equation}
with 
\begin{equation*}
R(\theta) = r_0 + 0.1 \, \sin( \omega \theta), \quad r_0 = 0.5, \quad \omega = 8.
\end{equation*}
$\Gamma$ is the zero level of the level set function $\phi(x) = \sqrt{x_1^2+x_2^2} - R(\theta)$ with $\theta(x) = \arctan(x_1/x_2)$.
The inner domain is $\Omega_1 := \{x \in \Omega, \phi(x) < 0 \}$. A sketch of the geometrical configuration is displayed in Figure \ref{fig:flower} (left) together with the initial \emph{unstructured} mesh which consists of $230$ triangles. 
Starting from that mesh, local refinements around the interface are successively performed 12 times. The final mesh consists of roughly $1.5$ million triangles. We consider different polynomial degrees $k$, where we always use the same polynomial degree for the approximation of $\phi$ and the deformation $\Psi_h$. 

In this example $\phi$ is not a signed distance function, s.t. $\text{dist}(\Gamma_h,\Gamma) \neq \Vert \phi \Vert_{\infty,\Gamma_h}$. However, there holds $\Vert \nabla \phi \Vert \geq 1$, s.t. we have $\text{dist}(\Gamma_h,\Gamma) \leq \Vert \phi \Vert_{\infty,\Gamma_h}$. In Figure \ref{fig:flower} the convergence of $\Vert \phi \Vert_{\infty,\Gamma_h}$ is depicted which is an upper bound for the distance.

If the mesh size $h$ in the vicinity of the interface is sufficiently small the geometrical error decreases with the (optimal) order $k+1$. For coarse grids, where the resolution of the interface is not sufficient the transformation has to be limited in order to garantuee shape regularity. On the corresponding grids we can not expect to observe a high order convergence which we also observe on the levels $0-2$. 
\subsection{A three-dimensional example: The gyroid}\label{sec:gyroid}
\begin{figure}[h!]
  \begin{center}
\includegraphics[trim=0.25cm 0cm 0.25cm 0cm, clip=true, width=0.475\textwidth]{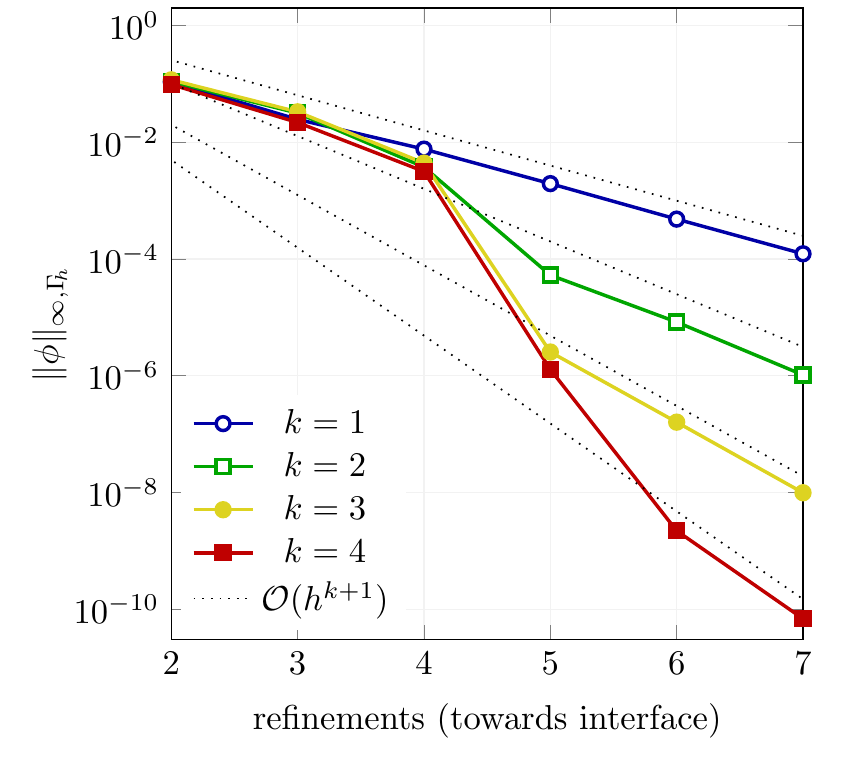} 
\hspace{0.025\textwidth}
  \includegraphics[trim=2.8cm 1cm 2.1cm 3cm, clip=true, width=0.425\textwidth]{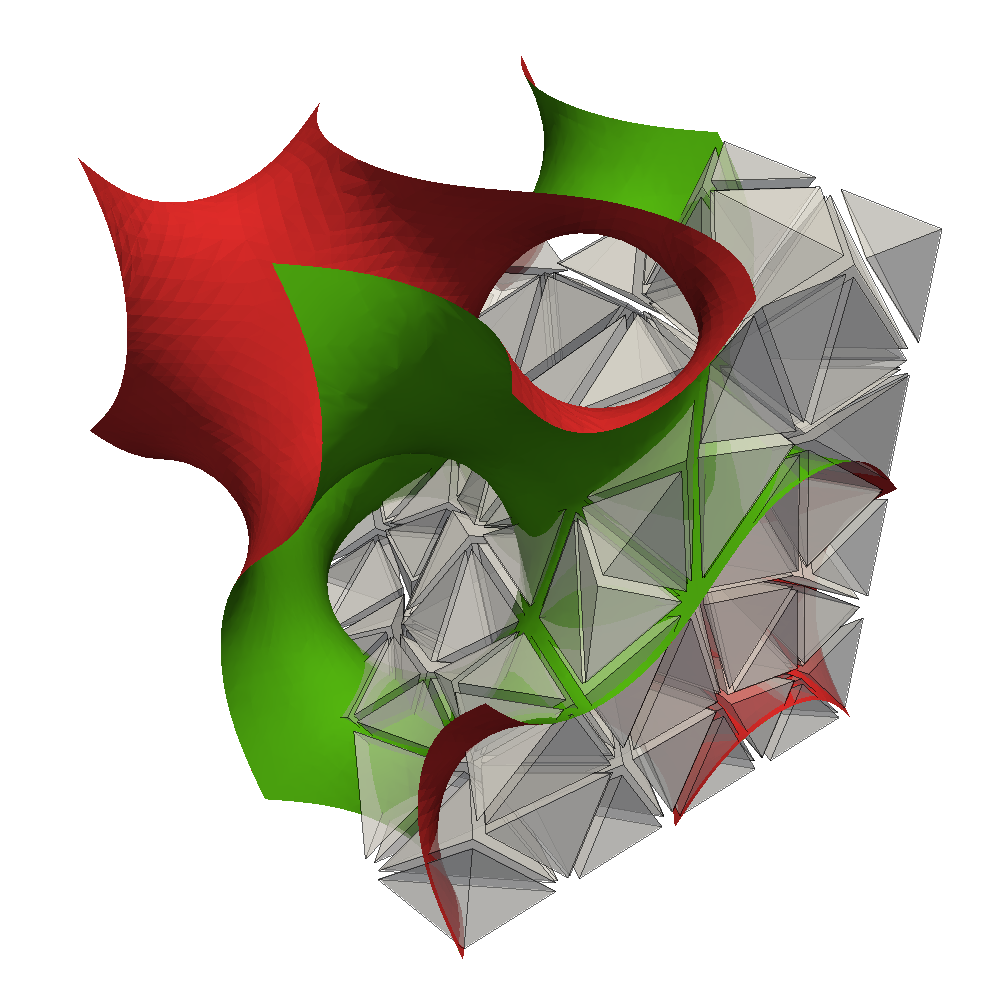}
\hspace{0.075\textwidth}
\end{center}
  \vspace*{-0.75cm}
  \caption{Mesh and convergence behavior of geometrical error for the gyroid example.}
\label{fig:gyroid}
\end{figure}

In this example we consider the cube $\Omega = [-1,1]^3$ with an interface prescribed by the zero level of the following level set function:
\begin{equation}
\phi(x,y,z) = 
\cos(\pi x_1) \sin(\pi x_2) + \cos(\pi x_2) \sin(\pi x_3) + \cos(\pi x_3) \sin(\pi x_1)
\end{equation}
A sketch of the interface and the initial \emph{unstructured} mesh, consisting of 577 tetrahedra, is depicted in Figure \ref{fig:gyroid}.
In a neighborhood of the interface $\Gamma$ there holds $\Vert \nabla \phi \Vert \geq 1$, s.t. we again have that the computed quantity $\Vert \phi \Vert_{\infty,\Gamma_h}$ is an upper bound for $\text{dist}(\Gamma_h,\Gamma)$. We consider seven adaptive refinements only towards the interface. The final mesh consists of around 40 million tetrahedra. Note that this time the interface is not contained in $\Omega$ such that $\Psi_h (\Omega) \neq \Omega$. However, for the approximation of $\Gamma$ and the evaluation of the error this is not important as $\phi$ is trivially extended to $\rr^3$.

The results in Figure \ref{fig:gyroid} reveal that the asymptotic behavior of the method is as good as in the previous two-dimensional example, the order of convergence $\mathcal{O}(h^{k+1})$ is clearly visible. However, the pre-asymptotic regime is significantly larger. This is due to the more difficult geometry. We observe that until refinement level $4$ the error drops only with $\mathcal{O}(h^{2})$. However, once the interface is resolved by the piecewise linear interface $\Gamma_1$, the limitation of the deformation is no longer necessary and the high order resolution of the interface kicks in. On the finest level, the accuracy between different polynomial degrees differs dramatically. An increase of the polynomial degree by one order results in 2 orders of magnitue decrease in the error. 

\subsection{Summary of numerical examples in this section}
These example represent configurations with realistic and challenging complexity. We conclude that the method is highly accurate if the interface is smooth and well resolved. 
If the interface is not well resolved, we at least obtain a robust method with the ``standard'' second order accuracy.

\section{Numerical examples II: A high order unfitted isoparametric finite element method for an interface problem} \label{sec:numexpde}
\subsection{Elliptic interface problem: The disk} \label{sec:disk}
\paragraph{Problem description:}
We consider the problem from \cite[section 2.5.1.4]{lehrenfeld2015diss}, where an elliptic interface problem of the following form is considered.
\begin{subequations} \label{eq:statNXFEM:modpoissoneq}
\begin{align}
\hspace*{2.0cm}  - \div (\alpha \nabla u) = & \, f 
 & \hspace*{-0.3cm} & \text{in}~~ \Omega_i, 
 & \hspace*{-0.3cm} & i=1,2, \hspace*{2.0cm}
\label{eq:statNXFEM:modpoissoneq1} \\
\jumpleft \alpha \nabla u \cdot n \jumpright   = & \, 0
 & \hspace*{-0.3cm} & \text{on}~~ \Gamma, 
 & \hspace*{-0.3cm} & 
 \label{eq:statNXFEM:modpoissoneq2} \\
\jumpleft \beta u \jumpright = & \, 0
 & \hspace*{-0.3cm} & \text{on}~~ \Gamma, 
 & \hspace*{-0.3cm} & 
 \label{eq:statNXFEM:modpoissoneq3} \\
 u = & \, g_D  
 & \hspace*{-0.3cm} & \text{on}~~\partial \Omega.
 & \hspace*{-0.3cm} & 
\label{eq:statNXFEM:modpoissoneq4} 
 \end{align}
\end{subequations}
Here $\alpha_i$, $\beta_i$, $i=1,2$ are domain-wise constants, $(\alpha_1, \alpha_2) = (2,1)$, $(\beta_1, \beta_2) = (1,3/2)$. We consider a circular interface $\Gamma$ with radius $R=0.6$ as interface $\Gamma$ in the square domain $[-1,1]^2$. The data $g_D$ and $f$ is chosen such that the solution to \eqref{eq:statNXFEM:modpoissoneq} is
\begin{align*}
&  u(x) = 
\left\{ 
\begin{array}{rr}
  \alpha_2 U(r(x)) + \beta_2, & x \in \Omega_1\\
  \alpha_1 U(r(x)) + \beta_1, & x \in \Omega_2\\
\end{array}
\right. \\
& \text{ with } \quad U(r) = \cos(\frac{\pi r^2}{2 R^2}) \text{ and } r(x) = \sqrt{x_1^2+x_2^2}
\end{align*}
\paragraph{The Nitsche-XFEM discretization:}
We consider an isoparametric version of the discretization in \cite[section 2]{lehrenfeld2015diss}, where a combination of an extended finite element (XFE) space combined with a Nitsche formulation for the unfitted interface is applied. The discretization is compactly described next. For a more thorough discussion of the underlying (planar) discretization we refer to \cite[chapter 2]{lehrenfeld2015diss}.

Let $\mathcal{R}_i$ be the restriction operator to domain $\Omega_{i,h}$ with $\Omega_{1,h}:= \{ I_h
\phi < 0 \}$ and $\Omega_{2,h}:= \{ I_h \phi > 0 \}$, such that $(\mathcal{R}_i u) |_{\Omega_j} =
\delta_{i,j} \cdot u|_{\Omega_i}$ with $\delta_{i,j}$ the Kronecker delta symbol. The extended
finite element space is then defined as $V_h^\Gamma = \mathcal{R}_1 V_h^k \oplus \mathcal{R}_2
V_h^k$. Basis functions from $V_h^\Gamma$ are continuous within the sub-domain $\Omega_{i,h}$ but
may be discontinuous across $\Gamma_1$. According to the transformation $\Psi_h$ we define the
isoparametric finite element space $\mathcal{V}_h^\Gamma = \{ v_h \circ \Psi_h^{-1} | v_h \in
V_h^\Gamma\}$ with correspondingly mapped shape functions. 

 The finite element space is not conforming w.r.t. the interface condition \eqref{eq:statNXFEM:modpoissoneq3}. To incorporate the interface condition we consider a Nitsche-type discretization: \\
Find $u_h \in \mathcal{V}_h^\Gamma$ such that 
\begin{equation}\label{eq:nitschexfem}
a(u_h,v_h) + N_h(u_h,v_h) = f(v_h) \text{ for all } v_h \in \mathcal{V}_h^\Gamma
\end{equation}
with $N_h(u,v) := N_h^c(u,v)+ N_h^c(v,u)+ N_h^s(u,v)$ and 
\begin{align*}
  a(u,v) &:= \sum_{i=1,2} \int_{\Psi_h(\Omega_{i,h})} \alpha_i \nabla u \nabla v \, dx, && f(v) := \int_{\Omega} f v \, dx \\ 
  N_h^c(u,v) &:= \int_{\Gamma_h} \averageleft - \alpha \nabla u \cdot \! n \averageright \jumpleft \beta v \jumpright \, ds, && 
  N_h^s(u,v) := \int_{\Gamma_h} \frac{ \bar\alpha \lambda k^2}{h} \jumpleft \beta u \jumpright \jumpleft \beta v \jumpright \, ds,
\end{align*}
with $\Gamma_h = \Psi_h(\Gamma_1)$, $\bar\alpha = \frac{1}{2}( \alpha_1 + \alpha_2 )$ and $\lambda$ a stabilization parameter that only depends on the shape regularity of the mesh.

The terms $a(u_h,v_h)$ and $f(v_h)$ ensure consistency inside the subdomains, $N^c_h(u_h,v_h)$ ensures consistency w.r.t. the interface condition \eqref{eq:statNXFEM:modpoissoneq2}. $N^c_h(v_h,u_h)$ render the discrete l.h.s. operator symmetric while keeping consistency and $N_s(u_h,v_h)$ ensures stability in a consistent way. 

An important detail in this discretization is the choice of the element-wise defined averaging operators $\averageleft \cdot \averageright$. We define the averaging weights $\kappa_i$, s.t. $\averageleft u \averageright := \kappa_1 u|_{\Omega_1} + \kappa_2 u|_{\Omega_2}$. We consider the ``Heaviside'' choice where $\kappa_i = 1$ if $|T_i| > \frac12 |T|$ and $\kappa_i = 0$ if $|T_i| < \frac12 |T|$. 
Here $|T_i|$ refers to the cut configuration on the undeformed mesh. This choice in the averaging renders the scheme in \eqref{eq:nitschexfem} stable for arbitrary polynomial degrees $k$, independent of the cut position of $\Gamma$. 
The interested reader is referred to the error analysis in \cite[section 2.3]{lehrenfeld2015diss}, specifically remark 2.3.1 therein. 

This test case has been considered for $k=1,2,3$ in \cite[section 2.5.1.4]{lehrenfeld2015diss} with
a piecewise planar approximation of the interface (on a eight times adaptively refined
subtriangulation). We now consider $k=1,2,3,4$ in $\mathcal{V}_h^\Gamma$ in combination with a level set function and a transformation $\Psi_h$ of the same polynomial degree. 

\paragraph{Numerical setup:}
We consider the mesh displayed in Figure \ref{fig:idea} as our initial mesh and successively apply uniform refinements of the mesh and measure the error in the $L^2$ norm, the broken $H^1$ (semi-) norm as well as the error in the interface condition \eqref{eq:statNXFEM:modpoissoneq3} in the $L^2(\Gamma_h)$-norm.
To solve the arising linear systems we use a sparse direct solver. We note that the arising linear systems are extremely ill-conditioned. We observed that in some cases, for $k>2$, the linear systems were so ill-conditioned that the initial residual of the linear system could only be reduced by a factor of $10^{-6}$. 
\begin{remark}
Efficient linear solvers for this kind of problem are rarely addressed in the literature, especially for high order unfitted finite elements. The only robust strategy (without additional stabilization) for linear systems arising from Nitsche-XFEM discretizations, that we know of, which is based on a rigorous analysis, is presented in \cite{reuskenlehrenfeld14}. In that paper however only the case $k=1$ is considered. The design and analysis of linear solvers for linear systems arising from Nitsche-XFEM discretizations remains a challenging task which requires further attention. 
\end{remark}

\paragraph{Numerical results:}
We note that no barrier step (cf. section \ref{sec:discretetrafo:shapereg}) has been necessary after one refinement. Hence, we expect to be in the asymptotic regime w.r.t. the geometry error such that $ \Vert \phi \Vert _{\infty,\Gamma_h} \leq \mathcal{O} (h^{k+1})$ holds. 
We observe the (optimal) convergence rates for the volume error that are predicted in the error
analysis in \cite[section 2.3]{lehrenfeld2015diss}. Note that error sources due to geometry
approximation and numerical integration have not been considered in that error analysis. The error
in the interface condition $\jumpleft \beta u \jumpright = 0$ converges with $\mathcal{O}(h^{k+1})$
which is even half an order better than predicted.

\begin{figure}[h!]
  \begin{center}
\hfill
\includegraphics[trim=0.25cm 0cm 0.25cm 0cm, clip=true, width=0.98\textwidth]{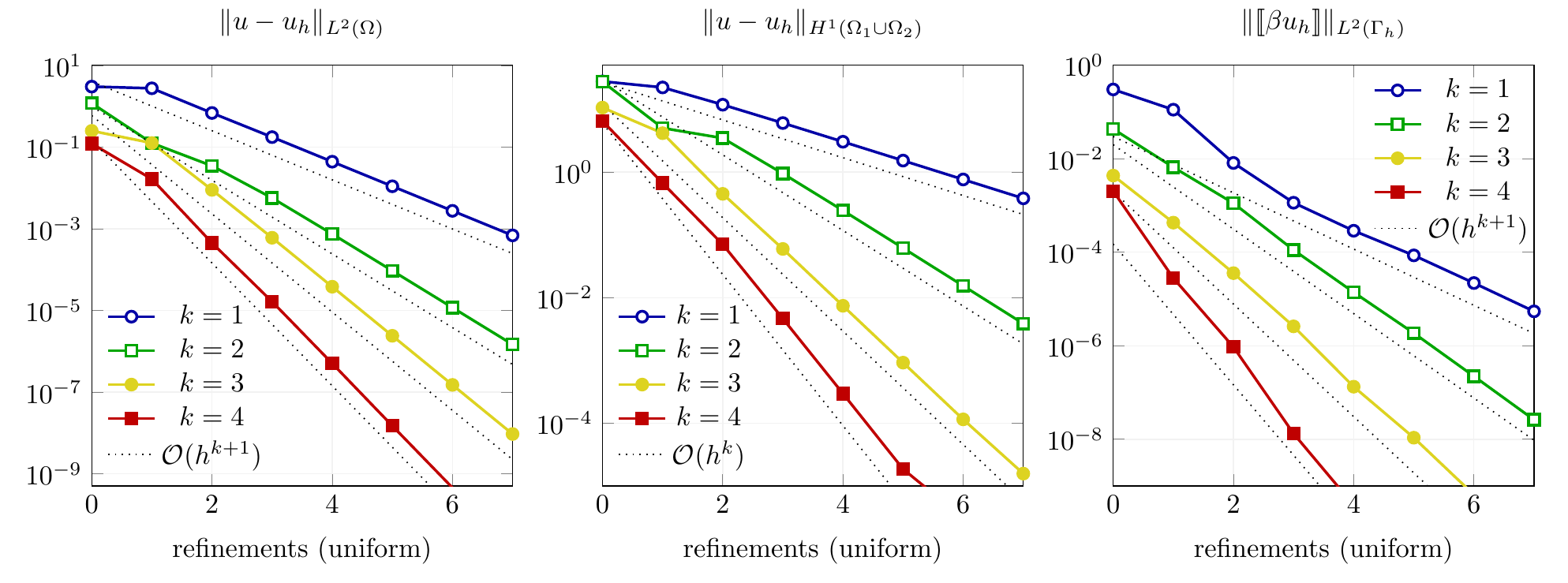} 
\end{center}
  \vspace*{-0.75cm}
  \caption{Error convergence in different norms for the circle example.}
\label{fig:circle}
\end{figure}

 \section*{Conclusion}
A new unfitted finite element method with a high order geometrical approximation of level set domains is presented and discussed in detail. The method is efficient and easy to implement. Numerical examples reveal its potential in handling complex geometries robust and highly accurate. 

As an example for the potential of the method, the method has been combined with a rather standard Nitsche-XFEM discretization for an unfitted interface problem and optimal order of convergence has been observed. The potential of the method goes beyond this specific example. As the method is geometry-based it allows to consider high order methods for a large range of unfitted finite element methods. However, many open questions remain. Stable discretizations and suitable linear solvers for high order unfitted finite element methods are difficult to develop and analyze. 

The presented method requires further attention. A rigorous error analysis of the presented method is missing, yet. We will address this in a forthcoming paper. Further, extensions to non-simplex meshes and time-dependent problems (with moving interfaces) are interesting topics for future studies.

\section*{Acknowledgements}
The author would like to express his appreciation to Arnold Reusken for his valuable and
constructive feedback on a first draft of this manuscript. The author also greatly appreciates the
fruitful discussion with Joachim Sch\"oberl on the topic of mesh transformations and shape
regularity in the context of this study.    

\bibliographystyle{elsarticle-num}
\bibliography{numint_lset}

\begin{thebibliography}{10}
\expandafter\ifx\csname url\endcsname\relax
  \def\url#1{\texttt{#1}}\fi
\expandafter\ifx\csname urlprefix\endcsname\relax\def\urlprefix{URL }\fi
\expandafter\ifx\csname href\endcsname\relax
  \def\href#1#2{#2} \def\path#1{#1}\fi

\bibitem{babuska73b}
I.~Babu\v{s}ka, The finite element method with penalty, Math. Comp. 27~(122)
  (1973) 221--228.

\bibitem{barrettelliott86}
J.~W. Barrett, C.~M. Elliott, Finite element approximation of the {Dirichlet}
  problem using the boundary penalty method, Numer. Math. 49~(4) (1986)
  343--366.

\bibitem{glowinskietal94}
R.~Glowinski, T.-W. Pan, J.~Periaux, A fictitious domain method for {Dirichlet}
  problem and applications, Comput. Meth. Appl. Mech. Eng. 111~(3--4) (1994)
  283--303.

\bibitem{burman2012fictitious}
E.~Burman, P.~Hansbo, Fictitious domain finite element methods using cut
  elements: {II}. a stabilized {Nitsche} method, Applied Numerical Mathematics
  62~(4) (2012) 328--341.

\bibitem{peskinmcqueen89}
C.~S. Peskin, D.~M. McQueen, A three-dimensional computational method for blood
  flow in the heart {I}. immersed elastic fibers in a viscous incompressible
  fluid, J. Comput. Phys. 81~(2) (1989) 372--405.

\bibitem{bastian2009unfitted}
P.~Bastian, C.~Engwer, An unfitted finite element method using discontinuous
  {Galerkin}, International journal for numerical methods in engineering
  79~(12) (2009) 1557--1576.

\bibitem{fries2010extended}
T.-P. Fries, T.~Belytschko, The extended/generalized finite element method: an
  overview of the method and its applications, International Journal for
  Numerical Methods in Engineering 84~(3) (2010) 253--304.

\bibitem{hansbo2002unfitted}
A.~Hansbo, P.~Hansbo, An unfitted finite element method, based on nitsche’s
  method, for elliptic interface problems, Computer methods in applied
  mechanics and engineering 191~(47) (2002) 5537--5552.

\bibitem{gross04}
S.~{Gro\ss}, V.~Reichelt, A.~Reusken, A finite element based level set method
  for two-phase incompressible flows, Comput. Visual. Sci. 9 (2006) 239--257.

\bibitem{massjung12}
R.~Massjung, An unfitted discontinuous {Galerkin} method applied to elliptic
  interface problems, SIAM J. Numer. Anal. 50~(6) (2012) 3134--3162.

\bibitem{Becker20093352}
R.~Becker, E.~Burman, P.~Hansbo, A {Nitsche} extended finite element method for
  incompressible elasticity with discontinuous modulus of elasticity, Computer
  Methods in Applied Mechanics and Engineering 198~(41–44) (2009) 3352 --
  3360.

\bibitem{olshanskii2009finite}
M.~A. Olshanskii, A.~Reusken, J.~Grande, A finite element method for elliptic
  equations on surfaces, SIAM journal on numerical analysis 47~(5) (2009)
  3339--3358.

\bibitem{lorensen1987marching}
W.~E. Lorensen, H.~E. Cline, Marching cubes: A high resolution 3d surface
  construction algorithm, in: ACM SIGGRAPH Computer Graphics, Vol.~21, ACM,
  1987, pp. 163--169.

\bibitem{naerland2014geometrychap5}
T.~A. N{\ae}rland, Geometry decomposition algorithms for the {Nitsche} method
  on unfitted geometries, Master's thesis, University of Oslo (2014).

\bibitem{mayer2009interface}
U.~M. Mayer, A.~Gerstenberger, W.~A. Wall, Interface handling for
  three-dimensional higher-order {XFEM}-computations in fluid--structure
  interaction, International Journal for Numerical Methods in Engineering
  79~(7) (2009) 846--869.

\bibitem{lehrenfeld2015nitsche}
C.~Lehrenfeld, The {Nitsche} {XFEM-DG} space-time method and its implementation
  in three space dimensions, SIAM Journal on Scientific Computing 37~(1) (2015)
  A245--A270.

\bibitem{lehrenfeld2015diss}
C.~Lehrenfeld, On a space-time extended finite element method for the solution
  of a class of two-phase mass transport problems, Ph.D. thesis, RWTH Aachen
  (February 2015).

\bibitem{DROPS}
S.~Gross, et~al., \href{http://www.igpm.rwth-aachen.de/DROPS}{{DROPS} package
  for simulation of two-phase flows} (2015).
\newline\urlprefix\url{http://www.igpm.rwth-aachen.de/DROPS}

\bibitem{engwer2012dune}
C.~Engwer, F.~Heimann, Dune-{UDG}: a cut-cell framework for unfitted
  discontinuous {Galerkin} methods, in: Advances in DUNE, Springer, 2012, pp.
  89--100.

\bibitem{burman2014cutfem}
E.~Burman, S.~Claus, P.~Hansbo, M.~G. Larson, A.~Massing,
  \href{http://dx.doi.org/10.1002/nme.4823}{{CutFEM}: Discretizing geometry and
  partial differential equations}, International Journal for Numerical Methods
  in Engineering.
\newline\urlprefix\url{http://dx.doi.org/10.1002/nme.4823}

\bibitem{renard2014getfem++}
Y.~Renard, J.~Pommier, \href{http://home.gna.org/getfem}{{GetFEM++}, an
  open-source finite element library} (2014).
\newline\urlprefix\url{http://home.gna.org/getfem}

\bibitem{carraro2015implementation}
T.~Carraro, S.~Wetterauer, On the implementation of the {eXtended} finite
  element method ({XFEM}) for interface problems, arXiv preprint
  arXiv:1507.04238.

\bibitem{chernyshenko2015adaptive}
A.~Y. Chernyshenko, M.~A. Olshanskii, An adaptive octree finite element method
  for pdes posed on surfaces, Computer Methods in Applied Mechanics and
  Engineering 291 (2015) 146--172.

\bibitem{muller2013highly}
B.~M{\"u}ller, F.~Kummer, M.~Oberlack, Highly accurate surface and volume
  integration on implicit domains by means of moment-fitting, International
  Journal for Numerical Methods in Engineering 96~(8) (2013) 512--528.

\bibitem{sudhakar2013quadrature}
Y.~Sudhakar, W.~A. Wall, Quadrature schemes for arbitrary convex/concave
  volumes and integration of weak form in enriched partition of unity methods,
  Computer Methods in Applied Mechanics and Engineering 258 (2013) 39--54.

\bibitem{saye2015hoquad}
R.~Saye, High-order quadrature method for implicitly defined surfaces and
  volumes in hyperrectangles, SIAM Journal on Scientific Computing 37~(2)
  (2015) A993--A1019.

\bibitem{burman2015cut}
E.~Burman, P.~Hansbo, M.~G. Larson, A cut finite element method with boundary
  value correction, arXiv preprint arXiv:1507.03096.

\bibitem{grande2014highorder}
J.~Grande, A.~Reusken, A higher order finite element method for partial
  differential equations on surfaces, Tech. Rep. 403, Institut f\"ur Geometrie
  und Praktische Mathematik, RWTH Aachen (2014).

\bibitem{cheng2010higher}
K.~W. Cheng, T.-P. Fries, Higher-order {XFEM} for curved strong and weak
  discontinuities, International Journal for Numerical Methods in Engineering
  82~(5) (2010) 564--590.

\bibitem{dreau2010studied}
K.~Dr{\'e}au, N.~Chevaugeon, N.~Mo{\"e}s, Studied {X-FEM} enrichment to handle
  material interfaces with higher order finite element, Computer Methods in
  Applied Mechanics and Engineering 199~(29) (2010) 1922--1936.

\bibitem{Basting2013228}
S.~Basting, M.~Weismann, A hybrid level set -- front tracking finite element
  approach for fluid--structure interaction and two-phase flow applications,
  Journal of Computational Physics 255 (2013) 228 -- 244.

\bibitem{oswald}
P.~Oswald, On a {BPX}-preconditioner for $\mathbb{P}_1$ elements, Computing 51
  (1993) 125--133.

\bibitem{schoeberl2014cpp11}
J.~Sch\"oberl, C++11 implementation of finite elements in {NGSolve}, Tech. Rep.
  ASC-2014-30, Institute for Analysis and Scientific Computing (September
  2014).

\bibitem{li95}
Z.~Li, A fast iterative algorithm for elliptic interface problems, SIAM Journal
  on Numerical Analysis 35 (1995) 230--254.

\bibitem{reuskenlehrenfeld14}
C.~Lehrenfeld, A.~Reusken, Optimal preconditioners for {Nitsche-XFEM}
  discretizations of interface problems, arXiv preprint arXiv:1408.2940.

\end{thebibliography}

\end{document}